\documentclass[12pt]{article}
\usepackage{amssymb}
\usepackage{latexsym}
\usepackage{amsmath}
\usepackage{amsthm}
\usepackage[OT2,T1]{fontenc}
\DeclareSymbolFont{cyrletters}{OT2}{wncyr}{m}{n}
\DeclareMathSymbol{\Sha}{\mathalpha}{cyrletters}{"58}
\newtheorem{theorem}{Theorem}[section]
\newtheorem{lemma}[theorem]{Lemma}
\newtheorem{proposition}[theorem]{Proposition}
\newtheorem{corollary}[theorem]{Corollary}
\newtheorem{definition}[theorem]{Definition}

\newtheorem{observation}[theorem]{Observation}
\newtheorem{question}[theorem]{Question}

\newtheorem{exercise}[theorem]{Exercise}

\newcommand{\lleg}[2]{\left(\!\!\left(\begin{array}{c} #1 \\ \hline #2 \end{array}\right)\!\!\right)}
\begin{document}
\title{Undecidability in number theory}
\author{Jochen Koenigsmann\\ Oxford}
\date{}
\maketitle
\newpage
\tableofcontents
\newpage
\section*{Introduction}
\addcontentsline{toc}{section}{Introduction}
These lectures are variations on a theme that is faintly echoed
in the following loosely connected counterpointing pairs:\\
\begin{center}
\begin{tabular}{rcl}
Euclid & versus & Diophantos\\
geometry & versus & arithmetic\\
decidability & versus & undecidability\\
Tarski & versus & G\"{o}del\\
Hilbert & versus & Matiyasevich
\end{tabular}
\end{center}
Let me explain how.

With his {\em Elements}
which in the Middle Ages
was the most popular `book' after the bible,
Euclid laid a foundation for modern mathematics already around 300 BC.
He introduced the axiomatic method
according to which every mathematical statement
has to be deduced from (very few) first principles (axioms)
that have to be so evident
that no further justification is required.
The paradigm for this is Euclidean geometry.

It wasn't quite so easy for arithmetic
(for good reasons as we know now).
In the 3rd century AD,
Diophantos of Alexandria,
often considered the greatest (if not only) algebraist of antique times,
tackled what we call today {\em diophantine equations},
that is, polynomial equations over the integers,
to be solved in integers.
Diophantos was the first to use symbols for unknowns,
for differences and for powers;
in short, he invented the polynomial.
He was the first to do arithmetic in its own right,
not just embedded into geometry (like, e.g., Pythagorean triples).
The goal was
to find a systematic method, a procedure, an algorithm
by which such diophantine equations could be solved
(like the well known formulas for quadratic equations).
One of the oldest, very efficient such algorithm
is the Euclidean (!) algorithm
for finding the greatest common divisor $\gcd(a,b)$
for two integers $a,b$.
This {\em is} a diophantine problem: for any intergers $a,b,c$,
$$c=\gcd(a,b)\Leftrightarrow\left\{
\begin{array}{l}
\mbox{the three diophantine equations}\\
cu=a,\; cv=b\mbox{ and }c=ax+by\\
\mbox{are solvable}
\end{array}\right.$$
(for even more surprising examples of mathematical problems
that are diophantine problems `in disguise', see section 3.4).

In his 10th problem from the famous list of 23 problems
presented to the Congress of Mathematicians in Paris in 1900,
David Hilbert, rather than asking for an algorithm
to produce solutions to diophantine equations,
asked for a more modest algorithm
that decides whether or not a given diophantine equation
{\em has} a solution.
In modern logic terminology,
such an algorithm would mean
that the existential 1st-order theory of $\mathbb{Z}$
(in the language of rings, ${\cal L}_{ring}:=\{ +,\times; 0,1\}$)
would be decidable (cf. section 3.1).

It would have been even more challenging
--- and quite in Hilbert's spirit ---
to show that the full 1st-order theory of $\mathbb{Z}$,
often simply called {\em arithmetic}, would be decidable.
However, in 1931, G\"{o}del showed
in the first of his two Incompleteness Theorems
that this is not the case:
{\em no algorithm can answer every arithmetic YES/NO-question correctly.}
It is called `Incompleteness Theorem' because it says
that every effectively (= algorithmically) producible
list of axioms true in $\mathbb{Z}$ is incomplete,
i.e., cannot axiomatize the full 1st-order theory of $\mathbb{Z}$
(section 1).

This undecidability result generalises to other number theoretic objects,
like all number fields (= finite extensions of $\mathbb{Q}$)
and their rings of integers,
by showing --- following Julia Robinson ---
that $\mathbb{Z}$ is 1st-order definable in any of these (section 2.3).
The key tools are the field $\mathbb{Q}_p$ of $p$-adic numbers (2.1)
and the Hasse-Minkowski Local-Gobal-Principle for quadratic forms (2.2).

In contrast, around the same time as G\"{o}del's Incompleteness Theorem,
Tarski showed that the other classical mathematical discipline,
geometry (at least elementary geometry), is decidable;
this holds true not only for Euclidean geometry,
but for all of algebraic geometry
where, when translated into cartesian coordinates,
geometric objects don't necessarily obey linear or quadratic equations,
but polynomial equations of arbitrary degree over $\mathbb{R}$ or $\mathbb{C}$:
The full 1st-order theory of $\mathbb{R}$
(and hence that of $\mathbb{C}$) is decidable.\footnote{Even though
Tarski may be better known for his decidability results than his undecidability results,
one should point out that his `Undefinability Theorem' (1936)
that {\em arithmetical truth cannot be defined in arithmetic} is very much in the spirit of G\"odel's Incompleteness Theorems
(in fact, it was discovered independently by G\"{o}del while proving these).
Tarski also proved undecidability of various other first-order theories, like, e.g., abstract projective geometry. This may put our very rough initial picture of the five counterpointing pairs
into a more accurate historical perspective.}

There is a whole zoo of interesting natural intermediate rings
between $\mathbb{Z}$ and $\mathbb{C}$
(or rather between $\mathbb{Z}$ and the field $\widetilde{\mathbb{Q}}$
of complex algebraic numbers).
To explore the boundaries within this zoo
of species that belong to the decidable world
(like $\widetilde{\mathbb{Q}}$,
or the field $\mathbb{Q}_p^{alg}:=\mathbb{Q}_p\cap\widetilde{\mathbb{Q}}$
of $p$-adic algebraic numbers or the field of totally real numbers
or the ring $\widetilde{\mathbb{Z}}$ of all algebraic integers)
and those on the undecidable side
(like number rings or number fields or the ring of totally real integers)
is a fascinating task 
with more open questions than answers (sections 2.4 and 2.5).

With the full 1st-order theory of $\mathbb{Z}$ being undecidable,
there still might be an algorithm to solve Hilbert's 10th Problem,
i.e., an effective decision procedure
for the existential 1st-order theory of $\mathbb{Z}$.
That this is also not the case is the celebrated result
due to Martin Davis, Hilary Putnam, Julia Robinson and Yuri Matiyasevich
(sometimes, for short, referred to as `Matiyasevich's Theorem'
as his contribution in 1970 was the last and perhaps most demanding):
{\em Hilbert's 10th Problem is unsolvable}
--- no algorithm can decide correctly
for all diophantine equations
whether or not they have integer solutions (section 3).

Maybe the most prominent open problem in the field
is the question whether Hilbert's 10th Problem can be solved over $\mathbb{Q}$,
i.e., whether there is an algorithm deciding solvability of
diophantine equations with solutions in $\mathbb{Q}$.
If we had an existential definition of $\mathbb{Z}$ in $\mathbb{Q}$
(which is still open) the answer would again be no
because then Hilbert's original 10th Problem over $\mathbb{Z}$
would be reducible to that over $\mathbb{Q}$,
and an algorithm for the latter would give one for the former,
contradicting Matiyasevich's Theorem.

Instead, in section 4, we reproduce the author's
{\em universal} definition of $\mathbb{Z}$ in $\mathbb{Q}$
which, at least in terms of logical complexity,
comes as close to the desired existential definition as one could get so far (4.1)
and, modulo either of two conjectures from arithmetic geometry,
as one ever possibly gets:
assuming Mazur's Conjecture or the Bombieri-Lang Conjecture,
there is no existential definition of $\mathbb{Z}$ in $\mathbb{Q}$ (4.3).

In section 5 we briefly discuss the question of full/existential decidability
for several other important rings,
not all arising from number theory.

In these notes we do not aim
at an encyclopedic survey of what has been achieved in the area,
nor do we provide full detailed proofs of the theorems treated
(each proof ought to be followed by an exercise: `fill in the gaps ...').
We rather try to point to the landmarks in the field
and their relative position,
to allow glimpses into the colourful variety of beautiful methods developed
for getting there.
Many (often, but not always long-standing) open problems are mentioned
to whet the appetite,
while the exercises provide working experience with some of the tools introduced.
What makes the topic really attractive, especially for graduate students,
is that most results don't use very heavy machinery,
are elementary in this sense, though, obviously,
people did have very good ideas.

I would like to express my warmest thanks
to Dugald Macpherson and Carlo Toffalori
for giving me the opportunity to hold these lectures
in the superb setting of Cetraro,
and to the enthusiastic audience
for their immense interest, their encouraging questions
and their critical remarks.
I am also very grateful to the anonymous referee and the editors
for their most valuable suggestions for improving on an earlier version
of these notes. 
\section{Decidability, Turing machines and G\"odel's 1st Incompleteness Theorem}
In this lecture we would like to sketch the proof of the first of G\"odel's celebrated two Incompleteness Theorems.
Denoting by ${\cal N}:=\langle\mathbb{N};+,\cdot;0,1\rangle$ the natural numbers
as ${\cal L}_{ring}$-structure, where ${\cal L}_{ring}:=\{+,\cdot;0,1\}$,
and by {\tt Th}$({\cal N})$ its 1st-order ${\cal L}_{ring}$-theory,
a weak version of the theorem is the following
\begin{theorem}[{[}G\"{o}d31{]}]
{\tt Th}$({\cal N})$ is undecidable.
\end{theorem}
I.e., there is {\em no} algorithm which,
on INPUT any ${\cal L}_{ring}$-sentence $\alpha$, gives
$$\mbox{OUTPUT }
\left\{\begin{array}{ll}
\mbox{\bf YES} & \mbox{if }\alpha\in\mbox{\tt Th}({\cal N}),\mbox{ i.e., }
{\cal N}\models\alpha\\
\mbox{\bf NO} & \mbox{otherwise}\end{array}\right.$$

In order to make this statement precise, we will define the notion of an {\em algorithm} using Turing machines.
There have been many alternative definitions
(via register machines, $\lambda$-calculus, recursive functions etc.)
all of which proved to be equivalent.
And, indeed, it is the credo of what has come to be called {\bf Church's Thesis}
that, no matter how we pin down an exact (and sensible) notion of algorithm,
it is going to be equivalent to the existing ones.
Whether or not one should take this as more than an empirical fact about the algorithms checked sofar, is an interesting philosophical question.
\subsection{Turing machines}
A {\em Turing machine} {\bf T} over a finite {\em alphabet}
$A =\{ a_1, a_2,\ldots ,a_n\}$
consists of a {\em tape}
$$\begin{array}{ccc|c|c|c|ccc} 
\hline
\cdot & \cdot & \cdot & a_i & a_j & a_k & \cdot & \cdot & \cdot\\   \hline
\end{array}$$
with infinitely many {\em cells}, each of which is either empty
(contains the {\em empty letter} $a_0$)
or contains exactly one $a_i$ ($1\leq i\leq n$).

In each step the {\em tape head} is on exactly one cell
and performs exactly one of the following four {\em operations}:
\begin{description}
\item[${\bf a_i}$]
type $a_i$ on the working cell ($0\leq i\leq n$)
\item[r] go to the next cell on the right
\item[l] go to the next cell on the left
\item[s] stop
\end{description}
The {\em program} (action table) of {\bf T}
is a finite sequence of lines of the shape\\
\begin{center}
{\bf z $a$ b z$^\prime$}
\end{center}
where $ a\in A\cup\{ a_0\}$,
where {\bf b} is one of the above operations,
and where {\bf z, z$^\prime$} are from a finite set
$Z=\{ z_1,\ldots ,z_m\}$ of {\em states}.

The program determines ${\bf T}$ by asking ${\bf T}$ to interpret
\begin{center}
{\bf z $a$ b z$^\prime$}
\end{center}
as `{\em if {\bf T} is in state {\bf z} and the tape head reads $a$
then do {\bf b} and go to state {\bf z$^\prime$}}'.

{\bf T} may stop on a given INPUT after finitely many steps
(and then the OUTPUT is what's on the tape then)
or it runs forever (with no OUTPUT given).
A decision algorithm always stops, by definition.
\subsection{Coding 1st-order ${\cal L}_{ring}$-formulas and Turing machines}
Let $A =\{ +,\cdot;0,1;\doteq,(,),\neg,\to,\forall,v,^\prime\}$
be the finite alphabet for 1st-order arithmetic,
thinking of the variable $v_n$ as the string $v^{\prime\prime\ldots\prime}$
of length $1+n$,
and assign to the finitely many elements of the disjoint union
$$A\cup Z\cup \{a_0,{\bf r,l,s}\}$$
distinct positive integers (their {\em codes}).
Code {\em formulas} via unique prime decomposition, e.g.,
if $\neg$, $0$, $\doteq$, $1$ have codes 2, 4, 1, 3 resp.,
the formula $\rho = \neg 0\doteq 1$ has code
$$\lceil\rho\rceil = 2^2\cdot 3^4\cdot 5^1\cdot 7^3 = 555660$$
and can be recovered from it.

Similarly, one can define a unique code $\lceil${\bf T}$\rceil$ for
(the program of) each {\em Turing machine} {\bf T} by coding the sequence of lines in the program.
\subsection{Proof of G\"odel's 1st Incompleteness Theorem (sketch)}
Suppose, for the sake of finding a contradiction,
that there is a Turing machine {\bf T}
which decides for any ${\cal L}_{ring}$-sentence $\alpha$
whether ${\cal N}\models \alpha$ or ${\cal N}\not\models\alpha$.
Then there is an arithmetic function $f_{\bf T}:\mathbb{N}\to\mathbb{N}$
describing ${\bf T}$ such that for any ${\cal L}_{ring}$-sentence $\alpha$,
$$f(\lceil \alpha\rceil)=\left\{\begin{array}{ll}
1 & \mbox{if }{\cal N}\models\alpha\\
2 & \mbox{if }{\cal N}\not\models\alpha\end{array}\right.$$
Here we call a function $f:\mathbb{N}\to\mathbb{N}$ {\em arithmetic}
if there is an ${\cal L}_{ring}$-formula $\phi (v_0,v_1)$
such that for any $n,m\in\mathbb{N}$,
$$f(n)=m\Leftrightarrow {\cal N}\models\phi (n,m).$$
($\phi$ is then called a {\em defining} formula for $f$.)
That there is such a defining formula $\phi_{\bf T}$ for $f_{\bf T}$
comes from the 1st-order fashion in which Turing programs come along;
the logical connectives and quantifiers translate into arithmetic operations.

So for any ${\cal L}_{ring}$-sentence $\alpha$ we have
$${\cal N}\models\phi_{\bf T}(\lceil\alpha\rceil,1)\Longleftrightarrow {\cal N}\models\alpha,$$
which means that we can `talk' about the truth of a sentence about ${\cal N}$ {\em inside} ${\cal N}$, and so we are in a position to simulate the liar's paradox:
Define $g:\mathbb{N}\to\mathbb{N}$ such that for all ${\cal L}_{ring}$-formulas $\rho (v_0)$
$$g(\lceil \rho\rceil)=\left\{\begin{array}{ll}
1 & \mbox{if }{\cal N}\models\neg\rho(\lceil\rho\rceil)\\
2 & \mbox{if }{\cal N}\models\rho(\lceil\rho\rceil) \end{array}\right.$$
Then $g$ is arithmetic, say with defining formula $\psi$, so that
$${\cal N}\models\psi(\lceil\rho\rceil,1) \Leftrightarrow {\cal N}\models\neg\rho (\lceil\rho\rceil).$$
Applied to the formula $\rho_0(v_0):=\psi (v_0,1)$ this gives
$${\cal N}\models\rho_0(\lceil\rho_0\rceil) \Leftrightarrow
{\cal N}\models\neg\rho_0(\lceil\rho_0\rceil),$$
the contradiction we were looking after.\qed
\medskip

Since there is an effective algorithm (a Turing machine) listing the Peano axioms, and since {\tt Th}$({\cal N})$ is complete we get the following immediate
\begin{corollary}
The Peano axioms don't axiomatise all of {\tt Th}$({\cal N})$.
\end{corollary} 
Nevertheless, a great many theorems about ${\cal N}$ do follow from the Peano axioms, and there has been an exciting controversy launched by Angus Macintyre as to whether Fermat's Last Theorem belongs
(cf. the Appendix in [Mac11]).
\medskip

Another immediate consequence is the following
\begin{corollary}
{\tt Th}$(\langle\mathbb Z;+,\cdot;0,1\rangle )$ is undecidable.
\end{corollary}
{\em Proof:} Otherwise, as the natural numbers are exactly the sums of four squares of integers, {\tt Th}$({\cal N})$ would be decidable.\qed
\section{Undecidability of number rings and fields}
\subsection{The field $\mathbb{Q}_p$ of $p$-adic numbers}
The field of $p$-adic numbers was discovered, or rather invented,
by Kurt Hensel over 100 years ago and has ever since played a crucial role in number theory.
It is, like the field $\mathbb{R}$ of real numbers,
the completion of $\mathbb{Q}$,
though not w.r.t. the ordinary (real) absolute value,
but rather w.r.t. a `$p$-adic' analogue (for each prime $p$ a different one).
This allows to bring new ($p$-adic) analytic methods into number theory
and to reduce some problems about number fields (which are so-called `global' fields) to the `local' fields $\mathbb{R}$ and $\mathbb{Q}_p$.
This is of particular interest in our context because, as we will see,
number fields are undecidable, whereas the fields of real or $p$-adic numbers are decidable.
So whenever a number theoretic problem is reducible
to a problem about the local fields $\mathbb{R}$ and $\mathbb{Q}_p$
(one then says that the problem satisfies a {\em Local-Global-Principle})
the number theoretic problem becomes decidable as well.\footnote{The use of the terms `local' and `global'
which one more typically encounters in analysis or in algebraic geometry
hints at a deep analogy between number theory and algebraic geometry,
more specifically between number fields
(i.e. finite extensions of $\mathbb{Q}$, the global fields of characteristic $0$)
and algebraic function fields in one variable over finite fields
(i.e., finite extensions of the field $\mathbb{F}_p(t)$ of rational functions over $\mathbb{F}_p$,
the global fields of positive characteristic).
It is one of the big open problems in model theory
whether or not the positive characteristic analogue of $\mathbb{Q}_p$,
i.e., the local field $\mathbb{F}_p((t))$ of (formal) Laurent series over $\mathbb{F}_p$,
is decidable.}

Let us fix a rational prime $p$.
The {\em $p$-adic valuation} $v_p$ on $\mathbb{Q}$ is defined by the formula
$$v_p(p^r\cdot \frac{m}{n})=r\mbox{ for any }r,m,n\in\mathbb{Z}\mbox{ with }
p\not\;\mid m\cdot n\neq 0,$$
with $v_p(0):=\infty$.
It is easy to check that this is a well defined valuation
(for background in valuation theory cf. [Dri13] in this volume or [EP05]).
The corresponding valuation ring is
$$\mathbb{Z}_{(p)} =\left\{ q\in\mathbb{Q}\;\vline\; v_p(q)\geq 0\right\}
= \left\{\frac{a}{b}\;\vline\; a\in\mathbb{Z}, p\not\;\mid b\in\mathbb{Z}\setminus\{ 0\}\right\}$$
with maximal ideal
$$p\mathbb{Z}_{(p)}=\left\{q\in\mathbb{Q}\;\vline\; v_p(q)> 0\right\}
= \left\{\frac{a}{b}\;\vline\; p\mid a\in\mathbb{Z}, p\not\;\mid b\in\mathbb{Z}\setminus\{ 0\}\right\}.$$
$v_p$ induces the {\em $p$-adic norm} $\mid . \mid_p$ on $\mathbb{Q}$ given by
$$\mid q\mid_p:=p^{-v_p(q)}\mbox{ for }q\neq 0,$$
with $\mid 0\mid_p :=0$.
Observe that the sequence $p,p^2,p^3,\ldots$ converges to $0$ w.r.t. $\mid . \mid_p$.

We can now define the {\em field $\mathbb{Q}_p$ of $p$-adic numbers}
as the completion of $\mathbb{Q}$ w.r.t. $\mid .\mid_p$,
i.e., the field obtained by taking the quotient of the ring of ($p$-adic) Cauchy sequences by the maximal ideal of ($p$-adic) zero sequences ({\tt version 1}).
Equivalently ({\tt version 2}), one may define
$$\mathbb{Q}_p:=\left\{\alpha =\sum_{\nu = n}^\infty a_\nu p^\nu\;\vline\;
n\in\mathbb{Z},a_\nu\in\{0,1,\ldots ,p-1\}\right\}$$
as the ring of formal Laurent series in powers of $p$ with coefficients from
$0,1,\ldots ,p-1$,
where addition is componentwise modulo $p$ starting with the lowest non-zero terms and carrying over whenever $a_\nu + b_\nu \geq p$ (so if this happens for the first time at $\nu$ the $(\nu +1)$-th coefficient becomes $1 + a_{\nu +1}+ b_{\nu +1}$ modulo $p$ etc.) ---
like adding decimals, but from the left.
Similarly, multiplication is like multiplying polynomials in the `unknown' $p$ with coefficients modulo $p$, and again carrying over whenever necessary.
\begin{exercise}
Show that
$-1 =\sum_{\nu =0}^\infty (p-1)p^\nu$ and 
$\frac{1}{1-p} = 1+p+p^2+\ldots$.
\end{exercise}
Given this presentation of $p$-adic numbers,
one defines the $p$-adic valuation on $\mathbb{Q}_p$,
again denoted by $v_p$, for any non-zero $\alpha\in\mathbb{Q}_p$ as
$$v_p(\alpha):=\min \{\nu\mid a_\nu\neq 0\}.$$
This is a prolongation of the $p$-adic valuation on $\mathbb{Q}$;
its value group is, obviously, still $\mathbb{Z}$,
and its valuation ring is
$$\mathbb{Z}_p:={\cal O}_{v_p} =\left\{\alpha\in\mathbb{Q}_p\;\vline\; v_p(\alpha)\geq 0\right\}
=\left\{\sum_{\nu=0}^\infty a_\nu p^\nu\;\vline\; a_\nu\in \{0,1,\ldots ,p-1\}\right\},$$
the {\em ring of $p$-adic integers}
with maximal ideal
$$p\mathbb{Z}_p =\left\{\alpha\in\mathbb{Q}_p\;\vline\; v_p(\alpha)> 0\right\}
=\left\{\sum_{\nu=1}^\infty a_\nu p^\nu\;\vline\; a_\nu\in \{0,1,\ldots ,p-1\}\right\},$$
and with residue field
$$\mathbb{Z}_p/p\mathbb{Z}_p\cong \mathbb{Z}_{(p)}/p\mathbb{Z}_{(p)}\cong
\mathbb{Z}/p\mathbb{Z} =\mathbb{F}_p.$$
By definition, $\mathbb{Q}$ is dense in $\mathbb{Q}_p$ w.r.t. the $p$-adic topology induced by (the two) $v_p$, and $\mathbb{Z}$ is dense in $\mathbb{Z}_p$:
in fact, $\mathbb{Z}_p$ is the completion of $\mathbb{Z}$
(w.r.t. the norm induced by $\mid .\mid_p$ on $\mathbb{Z}$).

Yet another way to think about the ring of $p$-adic integers
({\tt version 3}) is to view it as an inverse limit
$$\mathbb{Z}_p =\lim_{\leftarrow}\mathbb{Z}/p^n\mathbb{Z}$$
of the rings $\mathbb{Z}/p^n\mathbb{Z}$ w.r.t. the canonical projections
$\mathbb{Z}/p^n\mathbb{Z}\to\mathbb{Z}/p^m\mathbb{Z}$ for $m\leq n$.
Note that, for $n>1$, the rings $\mathbb{Z}/p^n\mathbb{Z}$ have zero divisors
whereas the projective limit $\mathbb{Z}_p$ becomes an integral domain
(with $\mathbb{Q}_p$ as field of fractions).
\begin{exercise}
Prove the equivalence of {\tt versions 1, 2, 3}.
\end{exercise}
One of the key facts about $\mathbb{Z}_p$ is {\em Hensel's Lemma}
which uses the analytic tool of Newton approximation
to find a precise zero of a polynomial, given an approximate zero.
For $\alpha\in\mathbb{Z}_p$,
let us denote its image under the canonical residue map
$\mathbb{Z}_p\to\mathbb{F}_p$ by $\overline{\alpha}$.
Similarly, we will write $\overline{f}$ for the image
of the polynomial $f\in\mathbb{Z}_p[X]$ under the
coefficientwise extension of the residue map to
$\mathbb{Z}_p[X]\to\mathbb{F}_p[X]$. 
\begin{lemma}[Hensel's Lemma]
{\bf Simple zeros lift:}
Let $f\in\mathbb{Z}_p[X]$ be a monic polynomial
and assume $\alpha\in\mathbb{Z}_p$ is such that
$\overline{\alpha}$ is a simple zero of $\overline{f}$
(i.e., $\overline{f}(\overline{\alpha})=0\neq\overline{f^\prime} (\overline{\alpha})$).
Then there is some $\beta\in\mathbb{Z}_p$ with
$f(\beta )=0$ and $\overline{\beta}=\overline{\alpha}$.
\end{lemma} 
\begin{exercise}
The proof is an adaptation of the proof in van den Dries' contribution to this volume, section 2.2,
the details being left to the reader as an exercise.
\end{exercise}
\noindent{\bf Example}
{\em If $p>2$ then every 1-unit,
that is, every $x\in 1+p\mathbb{Z}_p$ is a square:}\\
consider the polynomial $f(X) = X^2-x$ and let $\alpha =1$;
these satisfy the assumptions of Hensel's Lemma, and so
$f$ has a zero $\beta$ (with $\overline{\beta} =1$);
hence $x=\beta^2$.
\medskip\\
{\em Similarly, if one denotes by $\zeta_n$ a primitive $n$-th root of unity,
then $\zeta_{p-1}\in\mathbb{Z}_p$:}\\
the polynomial $X^{p-1} -1$ has $(p-1)$ distinct linear factors over $\mathbb{F}_p$, and so, by Hensel's Lemma, the same holds in $\mathbb{Z}_p$.

Thus we can write the multiplicative group of $\mathbb{Q}_p$ as a direct product of three `natural' subgroups:
$$\mathbb{Q}_p^\times = p^{\mathbb{Z}}\cdot\langle \zeta_{p-1}\rangle\cdot (1+p\mathbb{Z}_p)$$
From this one immediately reads off that, for $p>2$,
there are precisely four square classes
(elements in $\mathbb{Q}_p^\times/(\mathbb{Q}_p^\times)^2$),
represented by
$$1,\, p,\, \zeta_{p-1}\mbox{ and }p\zeta_{p-1}.$$
As a consequence, one obtains the following well-known 1st-order ${\cal L}_{ring}$-definition of $\mathbb{Z}_p$ in $\mathbb{Q}_p$:
$$\mathbb{Z}_p =\{x\in\mathbb{Q}_p\mid \exists y\in\mathbb{Q}_p\mbox{ such that } 1+px^2 = y^2\}.$$
\begin{exercise}
Show that, for $p=2$, a similar definition works with squares replaced by cubes.
\end{exercise}
\begin{exercise}
\label{quaternion}
Show that, for $p$ a prime $\equiv 3\mod 4$,
$$\mathbb{Z}_p=\{t\in\mathbb{Q}_p\mid\exists x,y,z\in\mathbb{Q}_p\mbox{ such that }
2+pt^2=x^2+y^2-pz^2\},$$
and, if $p\equiv 1\mod 4$ and $q\in\mathbb{N}$ a quadratic non-residue $\mod p$,
$$\mathbb{Z}_p=\{t\in\mathbb{Q}_p\mid\exists x,y,z\in\mathbb{Q}_p\mbox{ such that }
2+pqt^2=x^2+qy^2-pz^2\}.$$
\end{exercise}
With the valuation ring, also the maximal ideal, the residue field, the group of units, the value group and the valuation map all become interpretable in ${\cal L}_{ring}$. Thus, the axiomatization given below can be phrased entirely in ${\cal L}_{ring}$-terms.
Now here is a milestone in the model theory of $\mathbb{Q}_p$: 
\begin{theorem}[Ax-Kochen/Ershov]
{\tt Th}$(\mathbb{Q}_p)$ is decidable. It is effectively axiomatized by the following axioms:
\begin{itemize}
\item
$v_p$ is henselian
\item
the residue field of $v_p$ is $\mathbb{F}_p$
\item
the value group $\Gamma$ is a $\mathbb{Z}$-group,
i.e., $\Gamma\equiv \langle \mathbb{Z};+;0;<\rangle$
which can be axiomatized by saying that there is a minimal positive element
and that $[\Gamma :n\Gamma]=n$ for all $n$
\item
$v_p(p)$ is minimal positive
\end{itemize}
\end{theorem}
The proof, again, is similar to the proof of the other Ax-Kochen/Ershov Theorem presented in section 6 of van den Dries' contribution [Dri13] to this volume.

Fields elementarily equivalent to $\mathbb{Q}_p$ are called {\em $p$-adically closed}\footnote{Sometimes the term {\em $p$-adically closed} refers,
more generally, to fields elementarily equivalent to
{\em finite extensions} of $\mathbb{Q}_p$  --- cf. [PR84].}. 
\begin{exercise}
Check that fields which are relatively algebraically closed
in a $p$-adically closed field are again $p$-adically closed.
\end{exercise}
So, for example, the field
$$\mathbb{Q}_p^{alg}:=\mathbb{Q}_p\cap\widetilde{\mathbb{Q}}$$
of algebraic $p$-adic numbers is $p$-adically closed
(we use the notation $\widetilde{K}$ for the algebraic closure of $K$).
Note that $\mathbb{Q}_p^{alg}$ is countable while $\mathbb{Q}_p$ isn't.
\begin{exercise}
Show that $K:=\mathbb{Q}_p((\mathbb{Q}))$ (in the notation of [Dri13], after Definition 3.3,  this is
$\mathbb{Q}_p((t^{\mathbb{Q}}))$) is $p$-adically closed
and solve the mystery that,
on the one hand, $K$ and $\mathbb{Q}_p$ are elementarily equivalent,
on the other, they both have a henselian valuation
with the same residue field $\mathbb{Q}_p$,
but with non-elementarily equivalent value groups
($\mathbb{Q}$ for $K$ and $\{ 0\}$ for $\mathbb{Q}_p$).
\end{exercise}
\subsection{The Local-Global-Principle (LGP) for quadratic forms over $\mathbb{Q}$}
Let $q(X_1,\ldots ,X_n)$ be a quadratic form over $\mathbb{Q}$,
i.e., a homogeneous polynomial of degree $2$ in $\mathbb{Q}[X_1,\ldots ,X_n]$.
An element $a\in\mathbb{Q}$ is said to be {\em represented by $q$}
if there is $\overline{x}=(x_1,\ldots ,x_n)\in\mathbb{Q}^n$
such that $a=q(\overline{x})$.
\begin{theorem}[Hasse-Minkowski-Theorem]
\label{LGP}
A rational $a$ is represented by $q$ in $\mathbb{Q}$
if and only if
$a$ is represented by $q$ in all $\mathbb{Q}_p$ and in $\mathbb{R}$.
\end{theorem}
The proof is trivial for $n=1$,
it uses the so-called {\em geometry of numbers} for $n=2$,
it requires delicate case distinctions `modulo 8' for $n=3$ and $n=4$,
and then follows more easily by general quadratic form tricks for $n>4$.
(cf. [O'Me73], for an alternative proof using Hilbert symbols and quadratic reciprocity cf. [Ser73]).

The above LGP-principle is effective: by a simple linear transformation any quadratic form can be brought into diagonal form:
$$q(X_1,\ldots ,X_n)=a_1X_1^2 + \cdots + a_nX_n^2$$
Then the only primes $p$ where representability of $a$ by $q$ in $\mathbb{Q}_p$ need to be checked are $p=2$, $p=\infty$ (where $\mathbb{Q}_{\infty} :=\mathbb{R}$),
those $p$ where $v_p(a)\neq 0$, and those where $v_p(a_i)\neq 0$ for some $i\leq n$.
So only these finitely many primes need checking and, 
since all the $\mathbb{Q}_p$ and $\mathbb{R}$ are decidable,
the whole procedure is effective.

For {\em cubic} forms no such LGP holds:
By an example of Selmer ([Sel51]),
$5$ is represented by the cubic form
$$3X^3 + 4 Y^3$$
in $\mathbb{R}$ and in every $\mathbb{Q}_p$,
but not in $\mathbb{Q}$.
\subsection{Julia Robinson's definition for $\mathbb{Z}$ in $\mathbb{Q}$ and in other number fields}
Julia Robinson's contribution to questions of decidability in number theory is enormous.
Her first big result in this direction is
the 1st-order definability of $\mathbb{Z}$ (or $\mathbb{N}$) in $\mathbb{Q}$ 
from which the undecidability of {\tt Th}$(\mathbb{Q})$ immediately follows (1949), given G\"odel's 1st Incompleteness Theorem.
10 years later she extended this to arbitrary number fields.
Later she became heavily involved in Hilbert's 10th Problem (section 3).
The very appealing documentary
`Julia Robinson and Hilbert's Tenth Problem' by George Csicsery
came out in 2010 ([Csi10]).

To give Julia Robinson's explicit definition of $\mathbb{Z}$ in $\mathbb{Q}$,
let us introduce the following formulas:
for $a,b\in\mathbb{Q}^\times$ and $k\in\mathbb{Q}$, let
$$\phi (a,b,k) :=\exists x,y,z (2+abk^2+bz^2 = x^2 + ay^2)$$
and let, for $n\in\mathbb{Q}$,
$$\psi (n):=\forall a,b\neq 0
\left[ \{ \phi(a,b,0)\wedge\forall k\langle \phi(a,b,k)\to\phi(a,b,k+1)\rangle\}
\to\phi(a,b,n)\right]$$
\begin{theorem}[{[}Rob49{]}]
For any $n\in\mathbb{Q}$,
$$\mathbb{Q}\models\psi (n)\Leftrightarrow n\in \mathbb{Z}.$$
\end{theorem}
{\em Proof:}
The easy direction `$\Leftarrow$' follows, for $n\in\mathbb{N}$,
by the principle of induction,
and for $n\in\mathbb{Z}$, from the observation that
$\psi(n)\Leftrightarrow\psi(-n)$,
because $n$ occurs only squared in $\psi$.

For the non-trivial direction one first shows,
using Exercise \ref{quaternion} and Theorem \ref{LGP},
that, for a prime $p\equiv 3\mod 4$ and $k\in\mathbb{Q}$,
$$\phi(1,p,k)\Leftrightarrow v_p(k)\geq 0\mbox{ and }v_2(k)\geq 0,$$
and that, for primes $p,q$ with $p\equiv 1\mod 4$ and $q$ a quadratic non-residue $\mod p$,
$$\phi (q,p,k)\Leftrightarrow v_p(k)\geq 0\mbox{ and }v_q(k)\geq 0.$$
So, in either case, the $\{\ldots\}$-bit in $\psi$ is satisfied,
and, thus, for $\psi (n)$ to hold we must have
$\phi(1,p,n)$ for any prime $p\equiv 3\mod 4$
resp. $\phi(q,p,n)$ for any pair $p,q$ of primes in the second case.
But then, by the equivalences above, $v_p(n)\geq 0$ for any prime $p$,
and so $n\in\mathbb{Z}$.\qed
\begin{corollary}
{\tt Th}$(\mathbb{Q})$ is undecidable.
\end{corollary}
Let us recall, that number fields are finite extensions of $\mathbb{Q}$,
and that the {\em ring of integers in $K$},
denoted by ${\cal O}_K$, is the integral closure of $\mathbb{Z}$ in $K$,
i.e., the set of elements of $K$ satisfying a monic polynomial
with coefficients in $\mathbb{Z}$.
\begin{theorem}[{[}Rob59{]}]
\label{Rob59}
For any number field $K$,
${\cal O}_K$ is definable in $K$
and $\mathbb{Z}$ is definable in ${\cal O}_K$.
In particular, {\tt Th}$({\cal O}_K)$ and {\tt Th}$(K)$ are undecidable.
\end{theorem}
{\em Proof:}
The definition of ${\cal O}_K$ in $K$ proceeds along similar lines
as that of $\mathbb{Z}$ in $\mathbb{Q}$,
especially as the LGP for quadratic forms holds in arbitrary number fields.

That $\mathbb{N}$ is definable in ${\cal O}_K$ uses the fact
that for all non-zero $f\in {\cal O}_K$
there are only finitely many $a\in {\cal O}_K$ such that
$$a+1\mid f\wedge \ldots \wedge a+l\mid f,$$
where $l= [K:\mathbb{Q}]$.

Now define, for $a,f,g,h\in{\cal O}_K$,
$$\rho (a,f,g,h)\;:=\;f\not\doteq 0\wedge (a+1\mid f\wedge\ldots\wedge a+l\mid f)\wedge 1+ag\mid h$$
Then, for any $n\in {\cal O}_K$,
$$n\in\mathbb{N}\Leftrightarrow
\exists f,g,h\;\left[\rho (0,f,g,h)\wedge\forall a\{\rho (a,f,g,h)\to 
\langle\rho(a+1,f,g,h)\vee a\doteq n\rangle\}\right]$$
To prove the easy direction `$\Leftarrow$',
assume $n$ satisfies the right hand side.
By the fact above,
there are only finitely many $a$ with $\rho (a,f,g,h)$.
The inductive form of the definition ensures
that $\rho (0,f,g,h)$, $\rho(1,f,g,h),\ldots$
terminating only for $a=n$.
Therefore, $n$ must be a natural number.

For the converse direction `$\Rightarrow$',
assume $n\in\mathbb{N}$.
It suffices to find $f,g,h\in {\cal O}_K$ such that
$$\rho(a,f,g,h)\leftrightarrow a=0\vee a=1\vee\ldots\vee a=n.$$
Put $f:= (n+l)!$ and let
$S:=\{ a\in{\cal O}_K\mid a+1\mid f\wedge\ldots\wedge a+l\mid f\}$.
Then, by the fact above, $S$ is finite
and we can find some $g\in\mathbb{N}$ large enough so that,
for any two distinct $a,b\in S$, $a-b\mid g$
and, for any non-zero $a\in S$, $1+ag\not\,\mid 1$.
Then, for any distinct $a,b\in S$, $1+ag$ and $1+bg$ are relatively prime
(if there is a prime ideal of ${\cal O}_K$
containing both $1+ag$ and $1+bg$ then it contains $(a-b)g$, hence $g^2$ and so $g$, but it cannot contain both $g$ and $1+ag$).

Now put $h=(1+g)(1+2g)\cdots (1+ng)$.
Then $\rho(a,f,g,h)$ is satisfied for $a=1,\ldots ,n$.
If, however, there is some other $a\in S$
then $1+ag$ is not a unit and is prime to $h$.
Therefore, $1+ag\not\,\mid h$ and $\rho(a,f,g,h)$ does not hold.
\qed
\subsection{Totally real numbers}
There are many infinite algebraic extensions of $\mathbb{Q}$
(sometimes misleadingly called `infinite number fields')
which are also known to be undecidable.
In fact, most of them are: there are only countably many
decision algorithms, but uncountably many non-isomorphic,
and hence, in this case, non-elementarily equivalent
algebraic extensions of $\mathbb{Q}$.
To give an explicit example,
let $A$ be an undecidable (= non-recursive, cf. section 3.2) subset of the set of all primes and let
$$K:=\mathbb{Q}(\{\sqrt{p}\mid p\in A\}).$$
Then (the 1st-order theory of) $K$, and hence also ${\cal O}_K$,
is undecidable:
otherwise $A=\{p\in\mathbb{N}\mid p\mbox{ is prime and }\exists x\in K\;p=x^2\}$
would be decidable.

While this example seems artificial,
there is a number of `natural' infinite algebraic extensions of $\mathbb{Q}$ 
for which it makes sense to ask about decidability
of the field or of its ring of integers.
We will treat some of these in this section and the next,
and we will list some nice open problems in section 5.3.
The field of totally real numbers is special in that
its ring of integers is undecidable whereas the field is decidable.

{\em The field of totally real numbers} is defined to be
the maximal Galois extension $T$ of $\mathbb{Q}$ inside $\mathbb{R}$.
$T$ is an infinite algebraic extension of $\mathbb{Q}$,
the intersection of all real closures of $\mathbb{Q}$
(inside a fixed algebraic closure $\widetilde{\mathbb{Q}}$).
$T$ can also be thought of as the compositum of all finite extensions
$F/\mathbb{Q}$ for which all embeddings $F\hookrightarrow\mathbb{C}$ are real.

As for finite extensions of $\mathbb{Q}$
one defines ${\cal O}_T$, {\em the ring of integers of $T$},
as the integral closure of $\mathbb{Z}$ in $T$.
\begin{theorem}[{[}Rob62{]}]
\label{O_T}
{\tt Th}$({\cal O}_T)$ is undecidable.
\end{theorem}
Let us separate the key ingredients of the proof in the two lemmas below.
\begin{lemma}
\label{N in R}
Let $R$ be an integral domain with $\mathbb{N}\subseteq R$.
Let ${\cal F}\subseteq\wp (R)$ be a family of subsets of $R$
which is {\em arithmetically defined} (or {\em uniformly parametrised}),
say, by an ${\cal L}_{ring}$-formula $\phi(x;y_1,\ldots ,y_k)$,
i.e., for any $F\subseteq R$,
$$F\in {\cal F}\Leftrightarrow\exists\overline{y}\in R^k\;
\forall x\in R[x\in F\leftrightarrow\phi(x;\overline{y})].$$
Assume that all $F\in {\cal F}$ are finite
and each initial segment $\{0,1,\ldots ,n\}$ of $\mathbb{N}$
is in ${\cal F}$.
Then $\mathbb{N}$ is definable in $R$.
\end{lemma}
{\em Proof:}
For any $n\in R$,
$$n\in\mathbb{N}\Leftrightarrow
\exists\overline{y}\in R^k\;
\left[ \phi(0,\overline{y})\wedge\forall x\{ \phi(x,\overline{y})\to
\langle \phi(x+1,\overline{y})\vee x=n\rangle\}\right].$$
\qed

For $a,b\in T$, we use the notation `$a\ll b$'
to indicate that $a<b$ for any ordering $<$ on $T$.
Note that this is expressible by an ${\cal L}_{ring}$-formula:
$$a\ll b \Longleftrightarrow a\neq b\wedge\exists x\; b=a+x^2$$
(totally positive elements are always sums of squares,
and, in $T$, every sum of squares is a square). 

\begin{lemma}
\label{Kronecker}
$$\min \{M\in\mathbb{R}\mid\exists\infty\mbox{-ly many }t\in{\cal O}_T
\mbox{ s.t. }0\ll t\ll M\} = 4$$
\end{lemma}
{\em Proof:}
That there are infinitely many $t\in {\cal O}_T$
with $0\ll t\ll 4$ is easy: 
for any $n>1$ and any $n$-th root $\zeta_n$ of unity,
$t_n:=2 + \zeta_n + \zeta_n^{-1}$ has this property.

Conversely, any $t$ with this property is one of these $t_n$:
This follows from Kronecker's 1857 Theorem ([Kro57])
that {\em algebraic integers all of whose conjugates have
absolute value $\leq 1$ are roots of unity}:
if $\alpha=\alpha_1,\alpha_2,\ldots ,\alpha_n$
are all the conjugates of such an algebraic integer $\alpha$,
then, for any $k$, the coefficients $a_{k,s}$ of the polynomial
$$f_k(X):=(X-\alpha_1^k)\cdots (X-\alpha_n^k) =X^n+a_{k,n-1}X^{n-1}+\ldots +a_{k,0}$$
satisfy $\mid a_{k,s}\mid\leq \left( \begin{array}{c} n\\ s\end{array}\right)$;
being integers as well, there can only be finitely many such $a_{k,s}$,
hence only finitely many such $f_k$,
and so $\alpha^k=\alpha^l$ for some $k<l$,
making $\alpha$ an $(l-k)$-th root of unity.

From this one obtains that {\em totally real integers
all of whose conjugates have absolute value $\leq 2$
are of the shape $\alpha + \alpha^{-1}$
for some root of unity $\alpha$:}
let $\beta\in {\cal O}_T$ be such an element
with conjugates $\beta =\beta_1, \ldots, \beta_n$
and let $$\alpha=\frac{\beta}{2}+\sqrt{\frac{\beta^2}{4}-1};$$
then $\alpha$ is an algebraic integer
with  $\alpha^2-\beta\alpha +1=0$
and any conjugate $\alpha^\prime$ of $\alpha$ (over $\mathbb{Q}$)
satisfies $\alpha^{\prime 2}-\beta_i\alpha^\prime +1=0$ for some $i$;
as $\mid\beta_i\mid\leq 2$ (and, in fact, w.l.o.g., $<2$)
the two roots of this equation are the two complex conjugates of $\alpha^\prime$,
hence $$\mid\alpha^\prime\mid^2 =\frac{\beta_i^2}{4}+1-\frac{\beta_i^2}{4}=1,$$
so, by Kronecker's Theorem, $\alpha$ is a root of unity
and $\beta=\alpha +\alpha^{-1}$ is of the indicated shape.

Now the Lemma follows easily.\qed
\medskip

{\em Proof of Theorem \ref{O_T}:}
The family ${\cal F}\subseteq \wp ({\cal O}_T)$
defined by
$$\phi(x;p,q)\Leftrightarrow 0\ll qx\ll p\wedge p\ll 4q$$
contains, by Lemma \ref{Kronecker},
only finite sets, but arbitrarily large ones.
Hence, by a variant of Lemma \ref{N in R},
$\mathbb{N}$ is definable in ${\cal O}_T$.\qed
\subsection{Large algebraic extensions of $\mathbb{Q}$
and geometric LGP's}
As in the previous section,
we denote the field of totally real numbers by $T$.
In the literature, it is also often denoted by $\mathbb{Q}^{tot-r}$
or $\mathbb{Q}^{t.r.}$.
\begin{theorem}
\label{tot-r}
$T$ is {\em pseudo-real-closed} (`PRC'),
i.e. $T$ satisfies the following geometric LGP:
for each (affine) algebraic variety $V/T$,
$$V(T)\neq\emptyset \Leftrightarrow V(R)\neq\emptyset\mbox{ for all real closures }R\mbox{ of }T.$$
\end{theorem}
The theorem was first proved by Moret-Bailly ([Mor89])
using heavy machinery from algebraic geometry,
with a more elementary proof given later by Green, Pop and Roquette ([GPR95]).

By explicitly describing the structure of the absolute Galois group $G_T$ of $T$,
that is, the Galois group of the algebraic closure
$\widetilde{T}=\widetilde{\mathbb{Q}}$ of $T$ over $T$,
Fried, V\"olklein and Haran showed in [FVH94], using the above theorem:
\begin{theorem}
\label{FVH}
{\tt Th}$(T)$ is decidable.
\end{theorem} 
The axiomatization expresses the PRC property in elementary terms
(it is not at all obvious how to do this,
but it had long been established, e.g., in [Pre81])
as well as the fact that $G_T$ is the free (profinite) product of all $G_R$,
where $R$ runs through 
a set of representatives of the conjugacy classes of all real closures of $T$
(so each $G_R\cong\mathbb{Z}/2\mathbb{Z}$). 
The latter can be `axiomatized' via so-called embedding problems
in a similar fashion as, by a famous Theorem of Iwasawa,
free profinite groups of infinite countable rank can be characterized.

An interesting immediate consequence of this tension between
$T$ being decidable and ${\cal O}_T$ not, is the following:
\begin{corollary}
${\cal O}_T$ is not definable in $T$.
\end{corollary}

The decidability of (the 1st-order theory of) $T$ implies
that the field $T(\sqrt{-1})$ is decidable as well.
No answer is known to the following:
\begin{question}
Is ${\cal O}_{T(\sqrt{-1})}$ decidable?
\end{question}

A $p$-adic analogue of Theorem \ref{tot-r} and \ref{FVH}
was given by Pop in [Pop96]:
One defines the field $\mathbb{Q}^{tot-p}$,
the {\em field of totally $p$-adic numbers},
as the maximal Galois extension of $\mathbb{Q}$ inside $\mathbb{Q}_p$,
that is, the intersection of all conjugates of $\mathbb{Q}_p^{alg}$ over $\mathbb{Q}$.
Pop showed that $\mathbb{Q}^{tot-p}$ is {\em pseudo-$p$-adically-closed},
i.e. it satisfies an analogous LGP (with real closures replaced by $p$-adic closures), and that the absolute Galois group is similarly well behaved.
As a consequence, $\mathbb{Q}^{tot-p}$ is decidable.
We have no answer to the following
\begin{question}
Is ${\cal O}_{\mathbb{Q}^{tot-p}}$ decidable?
\end{question}

Let us close this section by mentioning another celebrated LGP,
that is, {\em Rumely's Local-Global-Principle} ([Rum86])
which concerns the ring $\widetilde{\mathbb{Z}}$
of all algebraic integers,
i.e., the integral closure of $\mathbb{Z}$
in the algebraic closure $\widetilde{\mathbb{Q}}$ of $\mathbb{Q}$:
\begin{theorem}
Let $V$ be an affine variety defined over $\widetilde{\mathbb{Z}}$.
Then
$$V(\widetilde{\mathbb{Z}})\neq\emptyset\Leftrightarrow
V({\cal O})\neq\emptyset
\mbox{ for all valuation rings }{\cal O}\mbox{ of }\widetilde{\mathbb{Q}}.$$ 
\end{theorem}
Using this, van den Dries ([Dri88]) showed the following Theorem via some quantifier elimination,
Prestel and Schmid ([PS90]) showed it via an explicit axiomatization:
\begin{theorem}
{\tt Th}$(\widetilde{\mathbb{Z}})$ is decidable.
\end{theorem}
Note that $\widetilde{\mathbb{Z}}$ is not definable in $\widetilde{\mathbb{Q}}$
(by quantifier elimination in ACF$_0$, every definable subset of
$\widetilde{\mathbb{Q}}$ is finite or cofinite),
so there is no cheap way of proving the above Theorem.

For a survey on geometric LGP's and many more results in this direction cf. [Dar00].
\section{Hilbert's 10th Problem and the DPRM-Theorem}
\subsection{The original problem and first generalisations}
In 1900, at the Conference of Mathematicians in Paris,
Hilbert presented his celebrated and influential list of 23 mathematical problems ([Hil00]).
One of them is
\medskip\\
{\bf Hilbert's 10th Problem (`H10')}
Find an algorithm which gives
on INPUT any $f(X_1, \ldots ,X_n)\in\mathbb{Z}[X_1,\ldots ,X_n]$
$$\mbox{OUTPUT} \left\{\begin{array}{ll}
\mbox{\bf YES} & \mbox{if }\exists\overline{x}\in\mathbb{Z}^n\mbox{ such that }f(\overline{x})=0\\
\mbox{\bf NO} & \mbox{else}
\end{array}\right.$$
Hilbert did not ask to prove that there is such an algorithm.
He was convinced that there should be one, and that it was all a question of producing it ---
one of those instances of Hilbert's optimism 
reflected in his famous slogan `{\em wir m\"{u}ssen wissen, wir werden wissen}'
(`{\em we must know, we will know}').
As it happens, Hilbert was too optimistic:
after previous work since the 50's by Martin Davis, Hilary Putnam and Julia Robinson,
in 1970, Yuri Matiyasevich showed  that there is no such algorithm
(Corollary \ref{noH10}).

The original formulation of Hilbert's 10th problem
was weaker than the standard version we have given above
in that he rather asked {\em `Given a polynomial $f$, find an algorithm ...}'.
So maybe one could have different algorithms depending on
the number of variables and the degree.
However, it is even possible to find a single polynomial
for which no such algorithm exists (Corollary \ref{strong})
--- this is essentially because there are universal Turing Machines.

One should, however, mention that, in the special case of $n=1$,
that is, for polynomials in one variable, there is an easy algorithm:
if, for some $x\in\mathbb{Z}$, $f(x)=0$ then $x\mid f(0)$;
hence one only has to check the finitely many divisors of $f(0)$.
Similarly, by the effective version of the Hasse-Minkowski-LGP (Theorem \ref{LGP}) and some extra integrality considerations,
one also has an algorithm for polynomials in an arbitrary number of variables, but of total degree $\leq 2$.
And, even if there is no general algorithm,
it is one of the major projects of computational arithmetic geometry
to exhibit other families of polynomials for which such algorithms exist. 

To conclude these introductory remarks
let us point in a different direction of generalizing Hilbert's 10th Problem,
namely, generalizing it to rings other than $\mathbb{Z}$:
If $R$ is an integral domain, there are two natural ways of generalizing {\bf H10}:
\medskip\\
\begin{center}
{\bf H10}/$R$ = {\bf H10} with the 2nd occurrence of $\mathbb{Z}$ replaced by $R$\\
{\bf H10}$^{+}/R$ = {\bf H10} with both occurrences of $\mathbb{Z}$ replaced by $R$
\end{center}
\begin{observation}
Let $R$ be an integral domain whose field of fractions does not contain the algebraic closure of the prime field ($\mathbb{F}_p$ resp. $\mathbb{Q}$). Then
$$\begin{array}{rcl}
\mbox{\bf H10}/R\mbox{ is solvable} & \Leftrightarrow & \mbox{\tt Th}_{\exists^+}(R)\mbox{ is decidable}\\
\mbox{\bf H10}^+/R\mbox{ is solvable} & \Leftrightarrow & \mbox{\tt Th}_{\exists^+}(\langle R;r\mid r\in R\rangle)\mbox{ is decidable,}
\end{array}$$
where {\tt Th}$_{\exists^+}$ denotes the positive existential theory
consisting of existential sentences where the quantifier-free part is
a conjunction of disjunctions of polynomial equations (no inequalities).
\end{observation}
Note that the language on the right hand side of the 2nd line contains a constant symbol for each $r\in R$.
\medskip\\
{\em Proof:}
`$\Leftarrow$' is obvious in both cases.
For `$\Rightarrow$' one has to see
that a disjunction of two polynomial equations
is equivalent to (another) single equation,
and, likewise, for conjunctions:
By our assumption we can find some monic $g\in\mathbb{Z}[X]$
of degree $>1$ which is irreducible over $R$.
Then, for any polynomials $f_1,f_2$ over $\mathbb{Z}$ resp. $R$
and for any tuple $\overline{x}$ over $R$,
$$\begin{array}{rcl}
f_1(\overline{x})=0\vee f_2(\overline{x})=0 & \Longleftrightarrow &
f_1(\overline{x})\cdot f_2(\overline{x}) =0\\
f_1(\overline{x})=0\wedge f_2(\overline{x})=0 & \Longleftrightarrow &
g(\frac{f_1(\overline{x})}{f_2(\overline{x})})\cdot f_2(\overline{x})^{\deg g} =0
\end{array}$$
\qed

Since in fields, inequalities can be expressed by a positive existential formula
($f(\overline{x})\neq 0\leftrightarrow\exists y\; f(\overline{x})\cdot y =1$),
we immediately obtain the following:
\begin{corollary}
Let $K$ be a field not containing the algebraic closure of the prime field. Then
$$\mbox{\bf H10}/K\mbox{ is solvable }\Leftrightarrow\mbox{\tt Th}_\exists(K)\mbox{ is decidable.}$$
\end{corollary}
In fact, the same is true for ${\cal O}_K$,
the ring of integers of a number field $K$:
\begin{exercise}
Show that, if $K$ is a number field,
$${\cal O}_K\models\forall x[x\neq 0\leftrightarrow\exists y\;x\mid (2y-1)(3y-1)].$$
Deduce that {\tt Th}$_\exists ({\cal O}_K) = $ {\tt Th}$_{\exists^+}({\cal O}_K)$.
\end{exercise}
One of the biggest open questions in the area is
\begin{question}
Is {\bf H10}$/\mathbb{Q}$ solvable?
\end{question} 
\subsection{Listable, recursive and diophantine sets}
A subset $A\subseteq\mathbb{Z}$ is called
\begin{itemize}
\item
{\em diophantine} if there is some $m\in\mathbb{N}$
and some polynomial $p\in\mathbb{Z}[T;X_1,\ldots ,X_m]$
such that
$$A=\{ a\in\mathbb{Z}\mid\exists\overline{x}\in\mathbb{Z}^m
\mbox{ with }p(a;\overline{x})=0\}$$
e.g., $\mathbb{N}$ is diophantine in $\mathbb{Z}$:
take $p=T-X_1^2-X_2^2-X_3^2-X_4^2$
\item
{\em listable} (= {\em recursively enumerable})
if there is an algorithm (= a Turing machine)
printing out the elements of $A$ (and only those),
e.g., the set $\mathbb{P}$ of primes or
$S:=\{a^3+b^3+c^3\mid a,b,c\in\mathbb{Z}\}$
\item
{\em recursive} (= {\em decidable})
if there is an algorithm deciding membership in $A$,
e.g., $\mathbb{N}$ and $\mathbb{P}$ are recursive,
about $S$ it is not known.
\end{itemize}
It is clear that every diophantine set is listable
and that every recursive set is listable.
That, conversely, every listable set is diophantine
is the content of the `DPRM-Theorem' (next section).

That not every listable set is recursive follows from the following
\begin{proposition}
[{\bf The Halting Problem of Computer Science is undecidable}]
There is no algorithm to decide whether a program (with code) $p$
halts on INPUT $x$.
\end{proposition}
{\em Proof:}
Otherwise define a new program $H$ by:
$$H\mbox{ halts on input }x\Leftrightarrow
x\mbox{ does not halt on input }x$$
(we identify $x$ with $\lceil x\rceil$).
For $x=H$ we are in trouble then.\qed

Using this, we find a listable, but non-decidable set:
$$A=\{ 2^p3^x\mid p\mbox{ halts on input }x\}$$
It is non-decidable by the proposition, but we can list it:
for $x,p\leq N$ print $2^p3^x$ if $p$ halts on input $x$
in $\leq N$ steps.
\subsection{The Davis-Putnam-Robinson-Matiyasevich\\
(= DPRM)--Theorem}
..., often for short referred to as Matiyasevich's Theorem,
is the following remarkable
\begin{theorem}
[{[}Mat70{]}, conjectured by Davis 1953, building on work of Davis, Putnam and Robinson]
\label{DPRM}
Every listable subset of $\mathbb{Z}$ is diophantine.
\end{theorem}
\begin{corollary}
\label{noH10}
Hilbert's 10th problem is unsolvable.
\end{corollary}
{\em Proof:}
The set $A$ at the end of the previous section is listable,
hence, by the Theorem, diophantine.
So there is some $m\in\mathbb{N}$
and some polynomial $p\in\mathbb{Z}[T;X_1,\ldots ,X_m]$
such that
$A=\{ a\in\mathbb{Z}\mid\exists\overline{x}\in\mathbb{Z}^m
\mbox{ with }p(a;\overline{x})=0\}$.
By construction, however, it is not decidable.
Hence there is no algorithm which decides, on input $t\in\mathbb{Z}$,
whether or not $p(t;\overline{x})=0$ has a solution $\overline{x}\in\mathbb{Z}^m$.
\qed
\medskip

We will only give a brief history and sketch of the proof of the Theorem.
For a full account of the history see the excellent survey article
[Mat00],
and, for a full self-contained proof (not always following the historic path),
see [Dav73].

Whether, in {\bf H10}, we ask for solutions in $\mathbb{Z}^n$
or in $\mathbb{N}^n$ doesn't make a difference:
as every integer is a difference of two natural numbers
and as every natural number is the sum of 4 squares of integers
we can easily transform a polynomial equation into another such that 
the former has integer solutions if and only if the latter has solutions
in natural numbers,
and we find a similar transformation for the other way round.
Following history
(and because $\mathbb{N}$ has the advantage of having a least element)
we will stick to finding solutions in $\mathbb{N}$.
\begin{theorem}
[{[}Dav53{]}]
\label{Davis}
If $A\subseteq\mathbb{N}$ is listable
then $A$ is {\em almost diophantine},
i.e., there is a polynomial $g\in\mathbb{Z}[T;\overline{X};Y,Z]$
such that for all $a\in \mathbb{N}$
$$a\in A\Leftrightarrow\exists z\forall y\leq z\;\exists\overline{x}\;
g(a;\overline{x};y,z)=0.$$
\end{theorem}
Note that every $A$ with such a presentation, later called `Davis normal form',
{\em is} listable.
\begin{exercise}
Observe that $2^{\mathbb{N}}$ is listable. Find an almost diophantine presentation.
\end{exercise}
A big challenge at the time was to find a diophantine presentation for $2^{\mathbb{N}}$, or, more generally, for exponentiation.
Julia Robinson showed that, in order to achieve this,
it suffices to find a diophantine relation `of exponential growth'
(what then went under the name `Julia Robinson-predicate'):
\begin{theorem}
[{[}Rob52{]}]
\label{Julia52}
There is a polynomial $q\in\mathbb{Z}[A,B,C;\overline{X}]$
such that for all $a,b,c\in\mathbb{N}$
$$a=b^c\Leftrightarrow\exists\overline{x}\; q(a,b,c;\overline{x})=0,$$
{\bf provided} there is a diophantine relation $J(u,v)$
{\em of exponential growth}, i.e., for all $u,v\in\mathbb{N}$
\begin{itemize}
\item
$J(u,v)\Rightarrow v< u^u$
\item
$\forall k\in\mathbb{N}\;\exists u,v$ with $J(u,v)$ and $v>u^k$.
\end{itemize}
\end{theorem}
That it is, indeed, enough to show that exponentiation is diophantine,
was then proved in the joint paper of Davis, Putnam and Julia Robinson\footnote{That we write the first name only for the lady is neither gallantry nor sexism:
the reason is that there are other Robinsons in the same area: Abraham Robinson, one of the founders of model theory, and the logician and number theorist Raphael Robinson, also Julia's husband}.
They used the term {\em exponential polynomial}
to refer to expressions obtained by applying the usual operations
of addition, multiplication and exponentiation
to integer coefficients and the variables:
\begin{theorem}
[{[}DPR61{]}]
\label{DPR}
If $A\subseteq\mathbb{N}$ is listable then there are exponential polynomials $E_L$ and $E_R$ such that, for all $a\in \mathbb{N}$,
$$a\in A\Leftrightarrow\exists\overline{x}\; E_L(a;\overline{x})=E_R(a;\overline{x}).$$
\end{theorem}
As predicted by Martin Davis, it then needed a young Russian mathematician,
Yuri Matiyasevich from St Petersborough, to fill the gap:
\begin{theorem}
[{[}Mat70{]}]
\label{Yuri}
There is a diophantine relation $J(u,v)$ of exponential growth.
\end{theorem}
The original proof used Fibonacci numbers.
In [D73], Davis gave a simpler proof using the so called {\em Pell equation}
$$x^2-dy^2=1,$$
where $d = a^2-1\in\mathbb{N}$ is a non-square.
It is not hard to verify that, for the number field $K:=\mathbb{Q}(\sqrt{d})$,
$$\begin{array}{rcl}
{\cal O}_K^\times & \supseteq & \{x+y\sqrt{d}\mid x,y\in\mathbb{Z}\mbox{ with }x^2-dy^2=1\}\\
 & = & \{\pm 1\}\cdot (a+\sqrt{a^2-1})^{\mathbb{Z}}
\end{array}$$
It then requires a series of elementary computations to actually check
that, on $\mathbb{N}\times\mathbb{N}$, the relation
$$J(x,y) \Leftrightarrow x^2-(a^2-1)y^2=1\Leftrightarrow x+y\sqrt{d}=(a+\sqrt{a^2-1})^m\mbox{ for some }m\in\mathbb{N}$$
is of exponential growth, or, at least, close to it.

\begin{center}
{\bf DPRM-Theorem \ref{DPRM} = Theorem \ref{Julia52} + \ref{DPR} + \ref{Yuri}}
\end{center}
\subsection{Consequences of the DPRM-Theorem}
One of the reasons why Davis' conjecture (the later DPRM-Theorem)
was considered dubious was the fact that it implies the existence of
a {\em prime producing polynomial}
($\mathbb{P}$ denotes the set of prime numbers):
\begin{corollary}
There is some $n\in\mathbb{N}$
and a polynomial $f\in\mathbb{Z}[X_1,\ldots ,X_n]$
such that $\mathbb{P} = f(\mathbb{Z}^n)\cap\mathbb{N}_{>0}$.
\end{corollary}
Today, even an explicit polynomial (with $n=10$) is known ([Mat81]).

{\em Proof:}
Obviously, $\mathbb{P}$ is listable.
Hence, by the DPRM-Theorem \ref{DPRM},
there is a polynomial $p\in\mathbb{Z}[T;\overline{X}]$ such that
$$\mathbb{P}=\{ a\in\mathbb{Z}\mid\exists\overline{x}\; p(a;\overline{x})=0\}.$$
But then the polynomial
$$f(T_1,\ldots ,T_4;\overline{X}):=[1-p(T_1^2+\ldots +T_4^2;\overline{X})](T_1^2 +\ldots +T_4^2)$$
does the job.\qed

As indicated at the beginning of this section,
Hilbert's 10th problem has a negative solution
even if one asks it for a single polynomial,
or, to put it less misleadingly,
for polynomials of a fixed shape,
in particular of a fixed number of variables and a fixed degree:
\begin{corollary}
\label{strong}
There is a polynomial $U\in\mathbb{Z}[T;\overline{X}]$
and an algorithm
producing, for each algorithm ${\cal A}$,
some $t_{\cal A}$ (a counterexample)
such that ${\cal A}$ fails to answer correctly whether
there is some $\overline{x}\in\mathbb{Z}^n$ with $U(t_{\cal A};\overline{x})=0$.
\end{corollary}
The proof relies on the fact that there are {\em universal}
recursive functions/ Turing machines which,
in addition to an INPUT-tuple $\overline{x}$,
take an INPUT-code $t_{\cal A}$,
and give as OUTPUT
the OUTPUT of the algorithm ${\cal A}$ on input $\overline{x}$.

Let us close this section by mentioning
that many famous mathematical problems
can be translated into diophantine problems.
For example, Goldbach's Conjecture that every even number $>2$
is the sum of two prime numbers:
Clearly, the set of counterexamples to this conjecture is listable,
hence, by DPRM, diophantine,
i.e., we can find a polynomial $g\in\mathbb{Z}[T;\overline{X}]$
such that, for any $t\in\mathbb{N}$,
$$g(t;\overline{x})=0$$
has a solution if and only if $t$ spoils the conjecture.
So Goldbach's Conjecture is equivalent to the statement that
$g(t; \overline{x})=0$ has no solution at all.

Similar translations can be found for Fermat's Last Theorem,
the Four Colour Theorem and the Riemann Hypothesis.
This may serve as an `explanation' why Hilbert's 10th problem
had to be unsolvable:
such difficult and diverse mathematical problems
cannot be expected to be solvable by just one universal process.
\section{Defining $\mathbb{Z}$ in $\mathbb{Q}$}
Hilbert's 10th problem over ${\mathbb Q}$,
i.e., the question whether {\tt Th}$_\exists ({\mathbb Q})$ is decidable, is still open.

If one had an {\em existential} (= {\em diophantine}) definition
of ${\mathbb Z}$ in ${\mathbb Q}$
(i.e., a definition by an existential 1st-order ${\cal L}_{ring}$-formula)
then  {\tt Th}$_\exists ({\mathbb Z})$ would be interpretable
in {\tt Th}$_\exists ({\mathbb Q})$,
and the answer would, by (for short) Matiyasevich's Theorem, again be no.
But it is still open
whether ${\mathbb Z}$ is existentially definable in ${\mathbb Q}$.

We have seen the earliest 1st-order definition of ${\mathbb Z}$ in ${\mathbb Q}$, due to Julia Robinson ([R49]), in section 2.3.
It can be expressed by an $\forall\exists\forall$-formula of the shape
$$\phi (t):\;  \forall x_1\forall x_2\exists y_1\ldots \exists y_7\forall z_1\ldots\forall z_6\; f(t;x_1,x_2;y_1,\ldots ,y_7; z_1,\ldots,z_6 )=0$$
for some $f\in {\mathbb Z}[t;x_1,x_2;y_1,\ldots ,y_7; z_1,\ldots,z_6]$,
i.e., for any $t\in {\mathbb Q}$,
$$t\in{\mathbb Z}\mbox{ iff } \phi (t) \mbox{ holds in }{\mathbb Q}.$$
In 2009, Bjorn Poonen ([P09a]) managed to find an $\forall\exists$-definition with 2 universal and 7 existential quantifiers
(earlier, in [CZ07], an $\forall\exists$-definition with just one universal quantifier was proved modulo an open conjecture on elliptic curves).
In this section we present our $\forall$-definition of ${\mathbb Z}$ in ${\mathbb Q}$:
\begin{theorem}
[{[}Koe10{]}]
\label{ZinQ}
There is a polynomial $g\in {\mathbb Z}[T;X_1,\ldots ,X_{418}]$ such that, for all $t\in {\mathbb Q}$,
$$t\in {\mathbb Z}\mbox{ iff }\forall \overline{x}\in {\mathbb Q}^{418}\; g(t;\overline{x})\neq 0.$$
\end{theorem}
If one measures logical complexity in terms of the number of changes of quantifiers then this is the simplest definition of ${\mathbb Z}$ in ${\mathbb Q}$,
and, in fact,  it is the simplest possible:
\begin{exercise}
Show that there is no quantifier-free definition of ${\mathbb Z}$ in ${\mathbb Q}$.
\end{exercise}

\begin{corollary}
${\mathbb Q}\setminus {\mathbb Z}$ is diophantine in ${\mathbb Q}$.
\end{corollary}
\begin{corollary}
$Th_{\forall\exists} ({\mathbb Q})$ is undecidable.
\end{corollary}
Theorem \ref{ZinQ} came somewhat unexpected because 
it does not give what one would like to have,
namely an existential definition of $\mathbb{Z}$ in $\mathbb{Q}$.
However, if one had the latter the former would follow:
\begin{observation}
If there is an existential definition of $\mathbb{Z}$ in $\mathbb{Q}$
then there is also a universal one.
\end{observation}
{\em Proof:}
If $\mathbb{Z}$ is diophantine in $\mathbb{Q}$ then so is
$$\mathbb{Q}\setminus\mathbb{Z}
=\{ x\in\mathbb{Q}\mid\exists m,n,a,b\in\mathbb{Z}\mbox{ with }
n\neq 0,\pm 1,\, am+bn=1\mbox{ and }m=xn\}$$\qed
\medskip

In fact, we will indicate in section 4.3 why we do not expect
there to be an existential definition of $\mathbb{Z}$ in $\mathbb{Q}$.
 
Using heavier machinery from number theory,
Jennifer Park has recently generalised Theorem \ref{ZinQ} to number fields:
\begin{theorem}
[{[}Par12{]}]
For any number field $K$, the ring of integers ${\cal O}_K$
is universally definable in $K$.
\end{theorem}
\subsection{Key steps in the proof of Theorem \ref{ZinQ}}
Like all previous definitions of ${\mathbb Z}$ in ${\mathbb Q}$,
we use the Hasse-Minkowski Local-Global-Principle for quadratic forms
(Theorem \ref{LGP}).
What is new in our approach is the use of the Quadratic Reciprocity Law and, inspired by the model theory of local fields, the transformation of some existential formulas into universal formulas.
\subsubsection*{{\em Step 1:} Poonen's diophantine definition of quaternionic semi-local rings}
The first step essentially copies Poonen's proof ([Poo09a]). We adopt his terminology:
\begin{definition}
For $a,b\in {\mathbb Q}^\times$, let
\begin{itemize}
\item
$H_{a,b}:=  {\mathbb Q}\cdot 1\oplus {\mathbb Q}\cdot \alpha\oplus {\mathbb Q}\cdot \beta\oplus {\mathbb Q}\cdot \alpha\beta$ be
the quaternion algebra over ${\mathbb Q}$ with multiplication defined by
$\alpha^2 =a$, $\beta^2 =b$ and $\alpha\beta = - \beta\alpha$,
\item
$\Delta_{a,b}:=\{l\in {\mathbb P}\cup\{\infty\}\mid H_{a,b}\otimes {\mathbb Q}_l\not\cong M_2({\mathbb Q}_l)\}$
the set of primes (including $\infty$) where $H_{a,b}$ does not split locally
(${\mathbb Q}_\infty :={\mathbb R}$)
--- $\Delta_{a,b}$ is always finite,
and $\Delta_{a,b}=\emptyset$ iff $a\in N(b)$,
i.e., $a$ is in the image of the norm map ${\mathbb Q} (\sqrt{b})\to {\mathbb Q}$,
\item
$S_{a,b}:= \{2x_1\in{\mathbb Q}\mid \exists x_2,x_3,x_4\in{\mathbb Q}:\,x_1^2-ax_2^2-bx_3^2+abx_4^2=1\}$ the set of traces of norm-$1$ elements of $H_{a,b}$, and
\item
$T_{a,b}:= S_{a,b}+S_{a,b}$ -- note that $T_{a,b}$ is an existentially defined subset of ${\mathbb Q}$.
\end{itemize}
\end{definition}\
\begin{lemma}
$T_{a,b}=\bigcap_{l\in\Delta_{a,b}}{\mathbb Z}_{(l)}$,
where, for $l\in \mathbb{P}$, $\mathbb{Z}_{(l)}=\mathbb{Z}_l\cap \mathbb{Q}$
is $\mathbb{Z}$ localised at $l$, and
${\mathbb Z}_{(\infty)}:=\{ x\in{\mathbb Q}\mid -4\leq x\leq 4\}$.
($T_{a,b}={\mathbb Q}$ if $\Delta_{a,b}=\emptyset$.) 
\end{lemma}
The proof follows essentially that of [Poo09a], Lemma 2.5, using Hensel's Lemma, the Hasse bound for the number of rational points on genus-$1$ curves over finite fields, and the local-global principle for quadratic forms.
Poonen then obtains his $\forall\exists$-definition of ${\mathbb Z}$ in ${\mathbb Q}$ from the fact that
$${\mathbb Z}=\bigcap_{l\in{\mathbb P}}{\mathbb Z}_{(l)}=\bigcap_{a,b>0} T_{a,b}.$$
Note that $\infty\not\in\Delta_{a,b}$ iff $a>0$ or $b>0$.
\subsubsection*{{\em Step 2:} Towards a uniform diophantine definition of all ${\mathbb Z}_{(p)}$'s in ${\mathbb Q}$}
We will present a diophantine definition
for the local rings ${\mathbb Z}_{(p)}={\mathbb Z}_p\cap {\mathbb Q}$
(i.e., $\mathbb{Z}$ localized at $p\mathbb{Z}$)
depending on the congruence of the prime $p$ modulo $8$,
and involving $p$ (and if $p\equiv 1\mod 8$ an auxiliary prime $q$) as a parameter.
However, since in any first-order definition of a subset of ${\mathbb Q}$
we can only quantify over the elements of ${\mathbb Q}$,
and not, e.g., over all primes,
we will allow arbitrary (non-zero) rationals $p$ and $q$
as parameters in the following definition. 
\begin{definition}
For $p,q\in {\mathbb Q}^\times$ let
\begin{itemize}
\item 
$R_p^{[3]} :=T_{-p,-p}+T_{2p,-p}$
\item 
$R_p^{[5]} :=T_{-2p,-p}+T_{2p,-p}$
\item 
$R_p^{[7]} :=T_{-p,-p}+T_{2p,p}$
\item
$R_{p,q}^{[1]} :=T_{2pq,q}+T_{-2pq,q}$
\end{itemize}
\end{definition}
The $R$'s are all existentially defined subrings of ${\mathbb Q}$ containing ${\mathbb Z}$,
since for any $a,b,c,d\in {\mathbb Q}^\times$
$$T_{a,b}+T_{c,d} = \bigcap_{l\in\Delta_{a,b}\cap\Delta_{c,d}} {\mathbb Z}_{(l)},$$
and since in each case at least one of $a,b,c,d$ is $>0$,
so $\infty\not\in\Delta_{a,b}\cap\Delta_{c,d}$.

\begin{definition}
\begin{enumerate}
\item[(a)]
${\mathbb P}^{[k]}:=\{ l\in {\mathbb P}\mid l\equiv k\mod 8\}$,
where $k=1,3,5$ or $7$
\item[(b)] 
For $p\in {\mathbb Q}^\times$, define
\begin{itemize}
\item 
${\mathbb P}(p):=\{ l\in {\mathbb P}\mid v_l(p)\mbox{ is odd}\}$,
where $v_l$ denotes the $l$-adic valuation on ${\mathbb Q}$
\item 
${\mathbb P}^{[k]}(p):={\mathbb P}(p)\cap {\mathbb P}^{[k]}$,
where $k=1,3,5$ or $7$
\item
$p\equiv_2 k\mod 8$ iff $p\in k+8{\mathbb Z}_{(2)}$,
where $k\in\{0,1,2,\ldots ,7\}$
\item
for $l$ a prime, the {\em generalized Legendre symbol} $\lleg{p}{l} =\pm 1$
to indicate whether or not the $l$-adic unit $pl^{-v_l(p)}$
is a square modulo $l$.
\end{itemize}
\end{enumerate}
\end{definition}
\begin{lemma}
\label{R_p}
\begin{enumerate}
\item[(a)] 
${\mathbb Z}_{(2)} = T_{3,3}+T_{2,5}$
\item [(b)]
For $p\in {\mathbb Q}^\times$ and $k=3,5$ or $7$,
if $p\equiv_2 k\mod 8$ then
$$R_p^{[k]}=\left\{\begin{array}{lll}
\bigcap_{l\in {\mathbb P}^{[k]}(p)} {\mathbb Z}_{(l)} & \mbox{if} & {\mathbb P}^{[k]}(p)\neq\emptyset\\
{\mathbb Q} & \mbox{if} & {\mathbb P}^{[k]}(p)=\emptyset
\end{array}\right.$$
In particular, if $p$ is a prime ($\equiv k\mod 8$)
then ${\mathbb Z}_{(p)} = R_p^{[k]}$.
\item [(c)]
For $p,q\in {\mathbb Q}^\times$ with $p\equiv_2 1\mod 8$ and $q\equiv_2 3\mod 8$,
$$R_{p,q}^{[1]}=\left\{\begin{array}{lll}
\bigcap_{l\in {\mathbb P}(p,q)} {\mathbb Z}_{(l)} & \mbox{if} & {\mathbb P}(p,q)\neq\emptyset\\
{\mathbb Q} & \mbox{if} & {\mathbb P}(p,q)=\emptyset
\end{array}\right.$$
where
$$l\in{\mathbb P}(p,q):\Leftrightarrow l\in\left\{\begin{array}{l}
{\mathbb P}(p)\setminus {\mathbb P}(q)\mbox{ with }\lleg{q}{l} =-1\mbox{, or}\\
{\mathbb P}(q)\setminus {\mathbb P}(p)\mbox{ with }\lleg{2p}{l} =\lleg{-2p}{l}=-1\mbox{, or}\\
{\mathbb P}(p)\cap {\mathbb P}(q)\mbox{ with }\lleg{2pq}{l} =\lleg{-2pq}{l}=-1\end{array}\right.$$
In particular, if $p$ is a prime $\equiv 1\mod 8$ and $q$ is a prime $\equiv 3\mod 8$ with $\left(\begin{array}{c}q\\ p\end{array}\right) =-1$ then ${\mathbb Z}_{(p)}=R_{p,q}^{[1]}$.
\end{enumerate}
\end{lemma}
\begin{corollary}
\label{ZviaRs} 
$${\mathbb Z} ={\mathbb Z}_{(2)}\cap\bigcap_{p,q\in{\mathbb Q}^\times}(R_p^{[3]}\cap R_p^{[5]}\cap R_p^{[7]}\cap R_{p,q}^{[1]})$$
\end{corollary}
The proof of the Lemma uses explicit norm computations for quadratic extensions of ${\mathbb Q}_2$, the Quadratic Reciprocity Law and the following
\begin{observation}
\label{linDelta}
For $a,b\in {\mathbb Q}^\times$ and for an odd prime $l$,
$$l\in\Delta_{a,b} \Leftrightarrow \left\{\begin{array}{l}
v_l(a)\mbox{ is odd, }v_l(b)\mbox{ is even, and }\lleg{b}{l} =-1\mbox{, or}\\
v_l(a)\mbox{ is even, }v_l(b)\mbox{ is odd, and }\lleg{a}{l} =-1\mbox{, or}\\
v_l(a)\mbox{ is odd, }v_l(b)\mbox{ is odd, and }\lleg{-ab}{l} =-1\\
\end{array}\right.$$
\end{observation}
\begin{corollary}
\label{dioph}
The following properties are diophantine properties for any $p\in {\mathbb Q}^\times$:
\begin{itemize}
\item 
$p\equiv_2 k\mod 8$ for $k\in\{0,1,2,\ldots ,7\}$ 
\item 
${\mathbb P}(p)\subseteq {\mathbb P}^{[1]}\cup {\mathbb P}^{[k]}$ for $k=3,5$ or $7$
\item
 ${\mathbb P}(p)\subseteq {\mathbb P}^{[1]}$
\end{itemize} 
\end{corollary}
\subsubsection*{{\em Step 3:} From existential to universal}
In Step 3, we try to find universal definitions for the $R$'s occurring in Corollary \ref{ZviaRs} imitating the local situation:
First one observes that for $R={\mathbb Z}_{(2)}$,
or for $R=R_p^{[k]}$ with $k=3, 5$ or 7, or for $R=R_{p,q}^{[1]}$,
the Jacobson radical $J(R)$ 
(which is defined as the intersection of all maximal ideals of $R$)
can be defined by an existential formula
using Observation \ref{linDelta}.
Now let
$$\widetilde{R}:=\{x\in {\mathbb Q}\mid\neg\exists y\in J(R)\mbox{ with }x\cdot y=1\}.$$
\begin{proposition}
\label{widetildeR}
\begin{enumerate}
\item[(a)]
$\widetilde{R}$ is defined by a {\em universal} formula in ${\mathbb Q}$.
\item[(b)]
If $R=\bigcap_{l\in{\mathbb P}\setminus R^\times}{\mathbb Z}_{(l)}$ then $\widetilde{R}=\bigcup_{l\in{\mathbb P}\setminus R^\times}{\mathbb Z}_{(l)}$,
provided ${\mathbb P}\setminus R^\times\neq\emptyset$,
i.e., provided $R\neq {\mathbb Q}$.
\item[(c)]
In particular, if $R={\mathbb Z}_{(l)}$ then $\widetilde{R}=R$.
\end{enumerate}
\end{proposition}
The proviso in (b), however, can be guaranteed
by diophantinely definable conditions on the parameters $p,q$:
\begin{lemma}
\label{RnotQ}
(a) Define for $k=1,3,5$ and $7$,
$$\begin{array}{ccl}
\Phi_k & := & \left\{ p\in {\mathbb Q}^\times\;\vline\; p\equiv_2 k\mod 8\mbox{ and } {\mathbb P}(p)\subseteq {\mathbb P}^{[1]}\cup {\mathbb P}^{[k]}\right\}
\vspace*{3mm}\\
\Psi & := & \left\{ (p,q)\in \Phi_1\times \Phi_3\;\vline\;
p\in 2\cdot ({\mathbb Q}^\times)^2\cdot (1+ J(R_q^{[3]}))\right\}.
\end{array}$$
Then $\Phi_k$ and $\Psi$ are diophantine in ${\mathbb Q}$.

(b) Assume that 
\begin{itemize}
\item
$R=R_p^{[k]}$ for $k= 3,5$ or $7$, where $p\in\Phi_k$, or
\item
$R=R_{p,q}^{[1]}$ where $(p,q)\in\Psi$.
\end{itemize}
Then $R\neq {\mathbb Q}$.
\end{lemma}
The proof of this lemma is somewhat involved,
though purely combinatorial,
playing with the Quadratic Reciprocity Law and Observation \ref{linDelta}.

The universal definition of ${\mathbb Z}$ in ${\mathbb Q}$ can now be read off the equation
$${\mathbb Z}=\widetilde{{\mathbb Z}_{(2)}}\cap(\bigcap_{k=3,5,7}\bigcap_{p\in\Phi_k}\widetilde{R_p^{[k]}})
\cap\bigcap_{(p,q)\in\Psi}\widetilde{R_{p,q}^{[1]}},$$
where $\Phi_k$ and $\Psi$ are the diophantine sets defined in Lemma \ref{RnotQ}.

The equation is valid by Lemma \ref{R_p}, Proposition \ref{widetildeR}(b), (c) and Lemma \ref{RnotQ}(b).
The definition is universal as one can see by spelling out the equation and applying Lemma \ref{widetildeR}(a) and Corollary \ref{dioph}: for any $t\in {\mathbb Q}$,
$$\begin{array}{lll}t\in{\mathbb Z} & \Leftrightarrow & t\in\widetilde{{\mathbb Z}_{(2)}}\wedge\\
 & &\forall p\bigwedge_{k=3,5,7} (t\in\widetilde{R_p^{[k]}}\vee p\not\in\Phi_k)\wedge\\
 & &\forall p,q (t\in \widetilde{R_{p,q}^{[1]}}\vee (p,q)\not\in\Psi)
\end{array}$$
Theorem \ref{ZinQ} is now obtained by diophantine routine arguments and counting quantifiers.
\subsection{More diophantine predicates in ${\mathbb Q}$}
From the results and techniques of section 4.1, one obtains new diophantine predicates in ${\mathbb Q}$. Among them are
\begin{itemize}
\item $x\not\in {\mathbb Q}^2$
\item $x\not\in N(y)$, where $N(y)$ is the image of the norm ${\mathbb Q}(\sqrt{y})\to{\mathbb Q}$
\end{itemize}
The first was also obtained in [Poo09b], using a deep result of Colliot-Th\'el\`ene et al. on Ch\^atelet surfaces --- our techniques are purely elementary.
\subsection{Why ${\mathbb Z}$ should not be diophantine in ${\mathbb Q}$}
There are two conjectures in arithmetic geometry that imply
that $\mathbb{Z}$ is not diophantine in $\mathbb{Q}$,
Mazur's Conjecture and, what one may call the Bombieri-Lang Conjecture.
\medskip\\
{\bf Mazur's Conjecture} ([Maz98])
{\em For any affine variety $V$ over $\mathbb{Q}$
the (real) topological closure of $V(\mathbb{Q})$ in $V(\mathbb{R})$
has only a finite number of connected components.}
\medskip\\
It is clear that, under this conjecture,
$\mathbb{Z}$ cannot be diophantine in $\mathbb{Q}$,
as the latter would mean
that $\mathbb{Z}$ is the projection of $V(\mathbb{Q})$
for some affine (not necessarily irreducible) variety $V$ over $\mathbb{Q}$,
but then, passing to the topological closure in $\mathbb{R}$,
$V(\mathbb{R})$ would have finitely many connected components
whereas the projection
(which is still the closed subset $\mathbb{Z}$ of $\mathbb{R}$)
has infinitely many - contradiction.
\medskip

The next conjecture,
though never explicitly formulated by Lang and Bombieri in this form,
may (arguably) be called `Bombieri-Lang Conjecture' (following [HS00]).
In order to state it we define,
given a projective algebraic variety $V$ over $\mathbb{Q}$,
the {\em special set} {\tt Sp}$(V)$
to be the Zariski closure of the union of all $\phi (A)$,
where $\phi:A\to V$ runs through all non-constant morphisms
from abelian varieties $A$ over $\mathbb{Q}$ to $V$.\medskip\\
{\bf Bombieri-Lang Conjecture}\\
{\em If $V$ is a projective variety over $\mathbb{Q}$
then $V(\mathbb{Q})\setminus(${\tt Sp}$(V)(\mathbb{Q}))$ is finite.}
\medskip\\
We shall use the following consequence of the conjecture:
\begin{lemma}
\label{Faltings}
Assume the Bombieri-Lang Conjecture.
Let $f\in\mathbb{Q}[x_1,\ldots ,x_{n+1}]\setminus\mathbb{Q}[x_1,\ldots ,x_n]$
be absolutely irreducible and
et $V=V(f)\subseteq\mathbb{A}^{n+1}$
be the affine hypersurface defined by $f$ over $\mathbb{Q}$.
Assume that $V(\mathbb{Q})$ is Zariski dense in $V$.
Let $\pi:\;\mathbb{A}^{n+1}\to\mathbb{A}^1$
be the projection on the 1st coordinate.
Then $V(\mathbb{Q})\cap\pi^{-1}(\mathbb{Q}\setminus\mathbb{Z})$
is also Zariski dense in $V$.
\end{lemma}
The proof uses a highly non-trivial finiteness result on integral points on abelian varieties by Faltings ([Fal91]).
\begin{theorem}
[{[}Koe10{]}]
\label{ZnotEinQ}
Assume the Bombieri-Lang Conjecture as stated above.
Then there is no infinite subset of ${\mathbb Z}$ existentially definable in ${\mathbb Q}$.
In particular, ${\mathbb Z}$ is not diophantine in ${\mathbb Q}$.
\end{theorem}
{\em Proof:}
Suppose $A\subseteq {\mathbb Z}$ is infinite and definable in ${\mathbb Q}$ by an existential formula $\phi_A(x)$ in the language of rings.
Replacing, if necessary, $A$ by $-A$,
we may assume that $A\cap\mathbb{N}$ is infinite.

Choose a countable proper elementary extension ${\mathbb Q}^\star$ of ${\mathbb Q}$ realizing the type $\{\phi_A(x)\wedge x>a\mid a\in A\}$
and let $A^\star =\{ x\in {\mathbb Q}^\star\mid\phi_A(x)\}$.
Then $A^\star$ contains some nonstandard natural number $x\in {\mathbb N}^\star\setminus {\mathbb N}$.
The map
{\small $\left\{\begin{array}{rcl} {\mathbb N} & \to & {\mathbb N}\\ n & \mapsto & 2^n\end{array}\right.$}
is definable in ${\mathbb N}$ and hence in ${\mathbb Q}$, so $2^x\in {\mathbb N}^\star$.
As $2^x$ is greater than any element algebraic over ${\mathbb Q}(x)$,
the elements $x, 2^x, 2^{2^x},\ldots $ are algebraically independent over ${\mathbb Q}$.
We therefore find an infinite countable transcendence base $\xi_1, \xi_2,\ldots $
of ${\mathbb Q}^\star$ over ${\mathbb Q}$ with $\xi_1\in A^\star$.

Let $K={\mathbb Q}(\xi_1, \xi_2, \ldots )$. 
As ${\mathbb Q}^\star$ is countable we find $\alpha_i\in {\mathbb Q}^\star$ ($i\in {\mathbb N}$) such that
$$K(\alpha_1)\subseteq K(\alpha_2)\subseteq\cdots \mbox{ with }\bigcup_{i=1}^\infty K(\alpha_i)={\mathbb Q}^\star,$$
where we may in addition assume that, for each $i\in {\mathbb N}$,
the minimal polynomial $f_i\in K[Z]$ of $\alpha_i$ over $K$
has coefficients in ${\mathbb Q}[\xi_1,\ldots ,\xi_i]$.
As ${\mathbb Q}$ is relatively algebraically closed in ${\mathbb Q}^\star$,
all the $f_i\in {\mathbb Q}[X_1,\ldots ,X_i,Z]$
are absolutely irreducible over ${\mathbb Q}$.

Now consider the following set of formulas in the free variables $x_1, x_2, \ldots$:
$$\begin{array}{rcl}
p = p(x_1,x_2,\ldots ) & := & \{ g(x_1,\ldots ,x_i)\neq 0\mid i\in {\mathbb N},
g\in {\mathbb Q} [x_1,\ldots ,x_i]\setminus\{ 0\}\}\\
 & & \cup \;\{\exists z\, f_i(x_1,\ldots ,x_i,z)=0\mid i\in {\mathbb N}\}\\
 & & \cup \;\{ x_1\mbox{ is not an integer}\}
\end{array}$$
Then $p$ is finitely realizable in ${\mathbb Q}$:
Let $p_0\subseteq p$ be finite and let $j$ be the highest index occurring in $p_0$ among the formulas from line 2.
Since the $K(\alpha_j)$ are linearly ordered by inclusion
all formulas from line 2 with index $<j$ follow from the one with index $j$.
Hence one only has to check that $V(f_j)$
has ${\mathbb Q}$-Zariski dense many ${\mathbb Q}$-rational points
$(x_1,\ldots ,x_j,z)\in {\mathbb A}^{j+1}$ with $x_1\not\in {\mathbb Z}$.
But this is, assuming the Bombieri-Lang Conjecture, exactly the conclusion of the above Lemma.
Note that $V(f_j)({\mathbb Q})$ is ${\mathbb Q}$-Zariski dense in $V(f_j)$ because there is a point
$(\xi_1,\ldots ,\xi_j,\alpha_j)\in V(f_j)({\mathbb Q}^\star )$
with $\xi_1,\ldots ,\xi_i$ algebraically independent over ${\mathbb Q}$.

Hence $p$ is a type that we can realize in some elementary extension
${\mathbb Q}^{\star\star}$ of ${\mathbb Q}$.
Calling the realizing $\omega$-tuple in ${\mathbb Q}^{\star\star}$ again
$\xi_1,\xi_2,\ldots $ our construction yields that we may view ${\mathbb Q}^\star$ as a subfield of ${\mathbb Q}^{\star\star}$.

But now $\xi_1\in A^\star\subseteq {\mathbb Z}^\star$ and $\xi_1\not\in {\mathbb Z}^{\star\star}$, hence $\xi_1\not\in A^{\star\star}$.
This implies that there is after all no existential definition for $A$ in ${\mathbb Q}$.\qed
\medskip

Let us conclude with a collection of closure properties for pairs of models of Th$({\mathbb Q})$ (in the ring language), one a substructure of the other,
which might have a bearing on the final (unconditional) answer
to the question whether or not ${\mathbb Z}$ is diophantine in ${\mathbb Q}$.
\begin{proposition}
Let ${\mathbb Q}^\star, {\mathbb Q}^{\star\star}$ be models of Th$({\mathbb Q})$ (i.e., elementary extensions of ${\mathbb Q}$) with ${\mathbb Q}^\star\subseteq {\mathbb Q}^{\star\star}$, and let ${\mathbb Z}^\star$ and ${\mathbb Z}^{\star\star}$ be their rings of integers. Then
\begin{enumerate}
\item[(a)]
${\mathbb Z}^{\star\star}\cap {\mathbb Q}^\star\subseteq {\mathbb Z}^\star$.
\item[(b)]
${\mathbb Z}^{\star\star}\cap {\mathbb Q}^\star$ is integrally closed in ${\mathbb Q}^\star$.
\item[(c)]
$({\mathbb Q}^{\star\star})^2\cap {\mathbb Q}^\star = ({\mathbb Q}^\star )^2$,
i.e. ${\mathbb Q}^\star$ is quadratically closed in ${\mathbb Q}^{\star\star}$.
\item[(d)]
If ${\mathbb Z}$ is diophantine in ${\mathbb Q}$
then ${\mathbb Z}^{\star\star}\cap {\mathbb Q}^\star = {\mathbb Z}^\star$
and ${\mathbb Q}^\star$ is algebraically closed in ${\mathbb Q}^{\star\star}$.
\item[(e)]
${\mathbb Q}$ is not model complete, i.e., 
there are $\mathbb{Q}^\star$ and $\mathbb{Q}^{\star\star}$ such that
${\mathbb Q}^\star$ is not existentially closed in ${\mathbb Q}^{\star\star}$.
\end{enumerate}
\end{proposition}
{\em Proof:}
(a) is an immediate consequence of our universal definition of ${\mathbb Z}$ in ${\mathbb Q}$.
The very same definition holds for ${\mathbb Z}^\star$ in ${\mathbb Q}^\star$
and for ${\mathbb Z}^{\star\star}$ in ${\mathbb Q}^{\star\star}$
(it is part of Th$({\mathbb Q})$ that all definitions of ${\mathbb Z}$ in ${\mathbb Q}$ are equivalent).
So if this universal formula holds for $x\in{\mathbb Z}^{\star\star}\cap{\mathbb Q}^\star$ in ${\mathbb Q}^{\star\star}$ it also holds in ${\mathbb Q}^\star$, i.e., $x\in {\mathbb Z}^\star$.

(b) is true because ${\mathbb Z}^{\star\star}$ is integrally closed in ${\mathbb Q}^{\star\star}$.

(c) follows from the fact that both being a square and, by section 4.2, not being a square are diophantine in ${\mathbb Q}$.

(d) If ${\mathbb Z}$ is diophantine in ${\mathbb Q}$
then ${\mathbb Z}^{\star\star}\cap {\mathbb Q}^\star\supseteq {\mathbb Z}^\star$ and hence equality holds, by (a).

To show that then also ${\mathbb Q}^\star$ is algebraically closed in ${\mathbb Q}^{\star\star}$, let us observe that, for each $n\in {\mathbb N}$,
$$A_n:=\{ (a_0,\ldots ,a_{n-1})\in{\mathbb Z}^n\mid\exists x\in {\mathbb Z}\mbox{ with }x^n+a_{n-1}x^{n-1} + \ldots + a_0=0\}$$
is decidable: zeros of polynomials in one variable are bounded in terms of their coefficients, so one only has to check finitely many $x\in {\mathbb Z}$.
In particular, by (for short) Matiyasevich's Theorem,
there is an $\exists$-formula $\phi (t_0,\ldots ,t_{n-1})$ such that
$${\mathbb Z}\models\forall t_0\ldots t_{n-1}\left(\{\forall x
[ x^n+t_{n-1}x^{n-1}+\ldots +t_0\neq 0]\}\leftrightarrow \phi (t_0,\ldots ,t_{n-1})\right).$$
Since both $A_n$ and its complement in ${\mathbb Z}^n$ are diophantine in ${\mathbb Z}$, the same holds in ${\mathbb Q}$, by our assumption of ${\mathbb Z}$ being diophantine in ${\mathbb Q}$,
i.e., $A_n^{\star\star}\cap ({\mathbb Q}^\star)^n = A_n^\star$.
As any finite extension of ${\mathbb Q}^\star$ is generated by an integral primitive element this implies that ${\mathbb Q}^\star$ is relatively algebraically closed in ${\mathbb Q}^{\star\star}$.

(e) Choose a recursively enumerable subset $A\subseteq {\mathbb Z}$ which is not decidable.
Then $B:={\mathbb Z}\setminus A$ is definable in ${\mathbb Z}$, and hence in ${\mathbb Q}$.
If $B$ were diophantine in ${\mathbb Q}$ it would be recursively enumerable. But then $A$ would be decidable: contradiction.

So not every definable subset of ${\mathbb Q}$ is diophantine in ${\mathbb Q}$,
and hence ${\mathbb Q}$ is not model complete.
Or, in other words, there are models ${\mathbb Q}^\star, {\mathbb Q}^{\star\star}$ of Th$({\mathbb Q})$ with ${\mathbb Q}^\star\subseteq {\mathbb Q}^{\star\star}$
where ${\mathbb Q}^\star$ is not existentially closed in ${\mathbb Q}^{\star\star}$.\qed
\medskip

We are confident that with similar methods as used in this paper one can show for an arbitrary prime $p$ that the unary predicate `$x\not\in {\mathbb Q}^p$' is also diophantine. This would imply that, in the setting of the Proposition, ${\mathbb Q}^\star$ is always radically closed in ${\mathbb Q}^{\star\star}$.
However, we have no bias towards an answer (let alone an answer) to the following (unconditional)
\begin{question}
For ${\mathbb Q}^\star\equiv {\mathbb Q}^{\star\star}\equiv {\mathbb Q}$ with
${\mathbb Q}^\star\subseteq {\mathbb Q}^{\star\star}$,
is ${\mathbb Q}^\star$ always algebraically closed in ${\mathbb Q}^{\star\star}$?
\end{question}
\section{Decidability and Hilbert's 10th Problem over other rings}
In this section we only report on major achievements under this heading
and on a small choice of big open problems.
There is a multitude of surveys on the subject, each with its own emphasis.
For the interested reader, let us mention at least some of them:
[RRo51], [Maz94], [Phe94], [PZ00], [Shl00], [Poo03], [Shl07] and [Poo08].
\subsection{Number rings}
For number rings and number fields,
the question of decidability has been answered in the negative
by Julia Robinson (Theorem \ref{Rob59}).
The question whether Hilbert's 10th Problem is solvable is much harder.
Given that we don't know the answer over $\mathbb{Q}$
(though almost everyone working in the field believes it to be no)
there is even less hope that we find the answer for arbitrary number fields in the near future.
For number {\em rings} the situation is much better.

Let $K$ be a number field with ring of integers ${\cal O}_K$.
Then Hilbert's 10th Problem could be shown to be unsolvable over ${\cal O}_K$
in the following cases:
\begin{itemize}
\item 
if $K$ is totally real (i.e., $K\subseteq T$)
or a quadratic extension of a totally real number field
([Den75], [DL78] and [Den80])
\item 
if $[K:\mathbb{Q}]\geq 3$ and $c_K=2$ ([Phe88])\footnote{$c_K$ denotes the {\em class number} of $K$, that is, the size of the ideal class group of $K$.
It measures how far ${\cal O}_K$ is from being a PID:
$c_K=1$ iff ${\cal O}_K$ is a PID, so $c_K=2$ is `the next best'.
It is not known whether there are infinitely many number fields with $c_K=1$.}.
\item 
if $K/\mathbb{Q}$ is abelian ([SS89]).
\end{itemize}
In each of the proofs the authors managed to find
an existential definition of $\mathbb{Z}$ in ${\cal O}_K$
using Pell-equations, the Hasse-Minkowski Local-Global Principle
(which holds in all number fields)
and ad hoc methods that are very specific to each of these special cases.

The hope for a uniform proof
of the existential undecidability of all number rings
only emerged when elliptic curves were brought into the game:
\begin{theorem}
[{[}Poo02{]}]
Let $K$ be a number field.
Assume\footnote{The set
$E(K)$ of $K$-rational points of $E$ is a finitely generated abelian group
isomorphic to the direct product of its torsion subgroup $E_{tor}(K)$
and a free abelian group of rank `{\tt rk}$(E(K))$'.}
there is an elliptic curve $E$ over $\mathbb{Q}$
with {\tt rk}$(E(\mathbb{Q}))= \mathbb{{\tt rk}}(E(K))=1$.
Then $\mathbb{Z}$ is existentially definable in ${\cal O}_K$
and so Hilbert's 10th Problem over ${\cal O}_K$ is unsolvable.
\end{theorem}
In his proof, Poonen uses divisibility relations
for denominators of $x$-coordinates of $n\cdot P$,
where $P\in E(K)\setminus E_{tor}(K)$
and $n\cdot P\in E(\mathbb{Q})$
(for a similar approach cf. [CPZ05]).

The assumption made in the theorem turns out to hold
modulo a generally believed conjecture,
the so called {\em Tate-Shafarevich Conjecture}.
For an elliptic curve $E$ over a number field $K$,
it refers to the {\em Tate-Shafarevich group}
(or {\em Shafarevich-Tate group})
$\Sha_{E/K}$, an abelian group defined via cohomology groups.
It measures the deviation from a local-global principle for
rational points on $E$.
\medskip\\
{\bf Tate-Shafarevich Conjecture}
{\em $\Sha_{E/K}$ is finite.}
\medskip\\
{\bf Weak Tate-Shafarevich Conjecture}
{\em $\dim_{\mathbb{F}_2}\Sha_{E/K}/2$ is even}
\medskip\\
The latter follows from the former due to the Cassels pairing
(Theorem 4.14 in [Sil86] which is an excellent reference on elliptic curves).
\begin{theorem}
[{[}MR10{]}]
Let $K$ be a number field.
Assume the weak Tate-Shafarevich Conjecture for all elliptic curves $E/K$.
Then there is an elliptic curve $E/\mathbb{Q}$
with {\tt rk}$(E(\mathbb{Q}))=\mbox{{\tt rk}}(E(K))=1$.
\end{theorem}
Taking those two theorems together one obtains immediately the following
\begin{corollary}
Let $K$ be a number field.
Assume the weak Tate-Shafarevich Conjecture for all elliptic curves $E/K$.
Then Hilbert's 10th Problem is unsolvable over ${\cal O}_K$.
\end{corollary}
\subsection{Function fields}
It is natural to ask decidability questions not only over number fields,
but also over global fields of positive characteristic,
i.e., algebraic function fields in one variable over finite fields,
and also, more generally, for function fields.

Hilbert's 10th Problem has been shown to be unsolvable for the following function fields:
\begin{itemize}
\item 
$\mathbb{R}(t)$ ([Den78])
\item 
$\mathbb{C}(t_1,t_2)$ ([KR92])
\item 
$\mathbb{F}_q(t)$ ([Phe91] and [Vid94])
\item
finite extensions of $\mathbb{F}_q(t)$ ([Shl92] and [Eis03])
\end{itemize}
The first two cases were achieved by existentially defining $\mathbb{Z}$
in the field, and then applying Matiyasevich's Theorem.
This is, clearly, not possible in the last two cases.
Instead of existentially {\em defining} $\mathbb{Z}$
the authors existentially {\em interpret} $\mathbb{Z}$
via elliptic curves:
the multiplication by $n$-map on an elliptic curve $E/K$
where $E(K)$ contains non-torsion points
easily gives a diophantine interpretation of the additive group
$\langle\mathbb{Z};+\rangle$.
The difficulty is to find an elliptic curve $E/K$
such that there is also an existential definition for multiplication
on that additive group.

For the ring of {\em polynomials} $\mathbb{F}_q[t]$,
Demeyer has even shown the analogue of the DPRM-Theorem:
listible subsets are diophantine ([Dem07]).

Generalizing earlier results ([Che84], [Dur86] and [Phe04]),
it is shown in [ES09], that the {\em full} first-order theory
of any function field of characteristic $>2$ is undecidable.

For analogues of Hilbert's 10th Problem for fields of meromorphic or analytic functions cf., e.g., [Rub95], [Vdau03] and [Pas13].
\subsection{Open problems}
Hilbert's 10th Problem is open for
\begin{itemize}
\item
$\mathbb{Q}$ and all number fields
\item
the ring of totally real integers ${\cal O}_T$ (cf. sections 2.4 and 2.5)
\end{itemize}
Hilbert's 10th Problem and full 1st-order decidability are open for
\begin{itemize}
\item 
$\mathbb{C}(t)$: this may well be considered the most annoying piece of our ignorance in the area.
On the other hand, $\mathbb{C}[t]$ and, in fact,
$R[t]$ for any integral domain $R$
is known to be existentially undecidable in ${\cal L}_{ring}\cup\{ t\}$
([Den78] and [Den84]).
\item 
$\mathbb{F}_p((t))$ and $\mathbb{F}_p[[t]]$ ---
in this case the answer to either question will be the same
for the field and the ring: in his recent thesis [Ans12],
Will Anscombe found a parameter-free existential definition
of $\mathbb{F}_p[[t]]$ in $\mathbb{F}_p((t))$.
In [DS03], Jan Denef and Hans Schoutens show
that Hilbert's 10th Problem {\em is} solvable,
if one assumes resolution of singularities in characteristic $p$. 
\item 
the field $\Omega$ of constructible numbers
(= the maximal pro-$2$ Galois extension of $\mathbb{Q}$).
What is known here is that ${\cal O}_\Omega$ is definable in $\Omega$ ([Vid99]),
and, more generally ([Vid00]), that for any prime $p$ and any pro-$p$ Galois extension $F$ of a number field, ${\cal O}_F$ is definable in $F$. As a consequence, the field of real numbers constructible with ruler and scale,
i.e., the maximal totally real Galois subextension $\Omega\cap T$ of $\Omega$ is undecidable.
\item
the maximal abelian extension $\mathbb{Q}^{ab}$ of $\mathbb{Q}$
and its ring of integers $\mathbb{Z}^{ab}$ ---
recall that, by the Kronecker-Weber Theorem,
$\mathbb{Q}^{ab}$ is the maximal cyclotomic extension of $\mathbb{Q}$,
obtained from $\mathbb{Q}$ by adjoining all roots of unity.
Here the answer may be related to the famous
Shafarevich Conjecture that the absolute Galois group
of $\mathbb{Q}^{ab}$ is a free profinite group
(if this is true decidability becomes more likely,
cf. the remarks following Theorem \ref{FVH}).

Let us remark that the ring of integers $\mathbb{Z}^{rab}$
of the field $\mathbb{Q}^{rab}:=\mathbb{Q}^{ab}\cap\mathbb{R}$
of real abelian algebraic numbers is undecidable,
by the identical proof as Theorem \ref{O_T}.
Note that $\mathbb{Q}^{ab}=\mathbb{Q}^{rab}(i)$
and that $\mathbb{Q}^{rab}$ is the fixed field of $\mathbb{Q}^{ab}$
under complex conjugation.
We do not know whether $\mathbb{Q}^{rab}$ is definable in $\mathbb{Q}^{ab}$,
but we conjecture that $\mathbb{Z}^{rab}$ is definable in $\mathbb{Z}^{ab}$
which would result in undecidability of $\mathbb{Z}^{ab}$.
\item
the maximal solvable extension $\mathbb{Q}^{solv}$ of $\mathbb{Q}$
and its ring of integers $\mathbb{Z}^{solv}$ ---
here the answer may be related to the longstanding open question
whether $\mathbb{Q}^{solv}$ is pseudo-algebraically-closed (PAC)
(Problem 10.16 in [FJ86], or Problem 11.5.9(a) in the 3rd edition;
a field $K$ is PAC if every absolutely irreducible algebraic variety
defined over $K$
has a $K$-rational point):
if the answer to this question is yes and if Shafarevich's Conjecture holds
then $\mathbb{Q}^{solv}$ {\em is} decidable,
axiomatized by being PAC, by the algebraic part and `by its absolute Galois group', the minimal normal subgroup of a free group with prosolvable quotient.
\item
the ring of integers of the field of totally $p$-adic numbers ---
here the question is whether there is an analogue of Kronecker's Theorem used in Lemma \ref{Kronecker}.
\end{itemize}
All known examples either have
both the full theory and the existential theory decidable
or both undecidable.
We have no answer to the following
\begin{question}
Is there a `naturally occurring' ring $R$
with {\tt Th}$_\exists(R)$ decidable, but
{\tt Th}$(R)$ not?
\end{question}
We are confident that one can build `unnatural' examples
using some Shelah-style construction
(though we haven't followed up on that).
A positive answer to an analogue of this question is Lipshitz's result that addition and divisibility over $\mathbb{Z}$ is $\exists$-decidable,
but $\forall\exists$-undecidable ([Lip78]).

Let us cconclude these notes by mentioning decidability questions
for fields that come up in other parts of this volume:
\begin{itemize}
\item 
Let $k$ be a field of characteristic $0$
and let $\Gamma$ be an ordered abelian group.
Then, by the Ax-Kochen/Ershov principle (see [Dri13]),
\begin{center}
\begin{tabular}{rcl}
$k((\Gamma))$ is decidable & $\Leftrightarrow$ & $k$ and $\Gamma$ are decidable\\
$k((\Gamma))$ is $\exists$-decidable & $\Leftrightarrow$ &
$k$ and $\Gamma$ are $\exists$-decidable
\end{tabular}
\end{center}
\item 
Let $\mathbb{R}_{exp}$ be the real exponential field
and let SC+ be the `souped up' version of Schanuel's Conjecture as in [MW96].
Then
\begin{center}
\begin{tabular}{rcl}
$\mathbb{R}_{exp}$ is decidable & $\Leftrightarrow$ & SC+ holds\\
 & $\Leftrightarrow$ & $\mathbb{R}_{exp}$ is $\exists$-decidable
\end{tabular}
\end{center}
\item 
The complex exponential field $\mathbb{C}_{exp}$ is, clearly, undecidable,
as $\mathbb{Z}$ is definable via the kernel of exponentiation
(this was already known to Tarski).
In fact, $\mathbb{Z}$ is even existentially definable in $\mathbb{C}_{exp}$:
In the 1980's, Angus Macintyre observed that, for any $x\in\mathbb{C}$,
$$x\in\mathbb{Q}\Leftrightarrow \exists t,v,u\left[ (v-u)t=1\wedge
e^v=e^u=1\wedge vx=u\right],$$
and in 2002, Mikl\'{o}s Laczkovich found, that for any $x\in\mathbb{C}$,
$$x\in\mathbb{Z}\Leftrightarrow x\in\mathbb{Q}\wedge\exists z\; (e^z=2\wedge e^{zx}\in\mathbb{Q}).$$
Hence, $\mathbb{C}_{exp}$ is even existentially undecidable
(cf. section 2.2 in [KMO12]).
\end{itemize} 

Mathematical Institute, 24-29 St Giles', Oxford OX1 3LB, UK\\
{\tt koenigsmann@maths.ox.ac.uk}
\end{document}